\documentclass[11pt,twoside,epsf]{article}
\usepackage{stmaryrd}
\usepackage{amsmath,latexsym,amssymb,amsfonts,amsbsy}
\usepackage{bm}
\usepackage{graphicx,subfigure}
\usepackage{cases}
\newfont{\bb}{msbm10}
\def\Bbb#1{\mbox{\bb #1}}

\def\tdiag{{\rm tridiag}}

\newtheorem{remark}{Remark}[section]
\newtheorem{theorem}{Theorem}[section]

\newtheorem{lemma}{Lemma}[section]

\newcommand{\reals}{\makebox{{\Bbb R}}}
\newcommand{\ceals}{\makebox{{\Bbb C}}}

\newcommand{\bbb}[1]{\text{\bf #1}}
\newcommand{\txt}[1]{\text{\rm #1}}

\makeatletter 

\@addtoreset{equation}{section}
\makeatother 

\begin{document}
\cleardoublepage
\pagestyle{myheadings}

\bibliographystyle{plain}

\title{The WR-HSS iteration method for a system of linear\\
differential equations and its applications to\\
the unsteady discrete elliptic problem \footnote{Supported by the
National Natural Science Foundation (No. 11101213 and No. 11401305),
P.R. China, and by the Natural Science Foundation of Jiangsu Province
(No. BK20141408), P.R. China.} }
\author
{Xi Yang
\footnote{Corresponding author at: Department of Mathematics,
Nanjing University of Aeronautics and Astronautics,
No. 29 Yudao Street, Nanjing 210016, P.R. China.
Email: yangxi@nuaa.edu.cn, yangxi@lsec.cc.ac.cn.}\\
{\it Department of Mathematics}\\
{\it Nanjing University of Aeronautics and Astronautics}\\
{\it No. 29 Yudao Street, Nanjing 210016, P.R. China}\\
}

\maketitle

\markboth{\small X. Yang}
{\small The WR-HSS method for discrete elliptic problem}

\begin{abstract}
We consider the numerical method for non-self-adjoint positive
definite linear differential equations, and its application to the
unsteady discrete elliptic problem, which is derived from spatial
discretization of the unsteady elliptic problem with Dirichlet
boundary condition. Based on the idea of the alternating direction
implicit (ADI) iteration technique and the Hermitian/skew-Hermitian
splitting (HSS), we establish a waveform relaxation (WR) iteration
method for solving the non-self-adjoint positive definite linear
differential equations, called the WR-HSS method. We analyze the
convergence property of the WR-HSS method, and prove that the WR-HSS
method is unconditionally convergent to the solution of the system
of linear differential equations. In addition, we derive the upper
bound of the contraction factor of the WR-HSS method in each iteration
which is only dependent on the Hermitian part of the corresponding
non-self-adjoint
positive definite linear differential operator. Finally, the applications
of the WR-HSS method to the unsteady discrete elliptic problem demonstrate
its effectiveness and the correctness of the theoretical results.

\bigskip


\end{abstract}

\section{Introduction}
\label{sec-introduction}

We consider the numerical solution of the following unsteady elliptic
problem (second-order parabolic equation),
\begin{eqnarray}
\label{uns-elliptic} \left\{
\begin{array}{l}
  \frac{\partial u(x,t)}{\partial t}-\nabla \cdot [a(x,t)\nabla u(x,t)] +
  \sum _{j=1}^{d} \frac{\partial}{\partial x_j}(p(x,t)u(x,t)) = q(x,t),\,x\in\Omega,\,t\in [0,\mbox{T}],\\
  u(x,0) = u_0(x),\,x\in\Omega,\\
  u(x,t)=v(x,t),\,x\in\partial\Omega,\,t\in [0,\mbox{T}]
\end{array}
\right.
\end{eqnarray}
with $\Omega$ being a plurirectangle of $\reals^d$, $\partial\Omega$ being
the boundary of the domain $\Omega$, $\mbox{T}$ (possibly infinite) being
the upper bound of the time interval, $a(x,t)$ being a uniformly
positive function and $p(x,t)$ denoting the Reynolds function. Specifically,
plurirectangle here means a connected union of rectangles in $d$-dimensions
with edges parallel to the axes. The above equation is important for various
reasons \cite{Elman-05}. As well as describing many significant physical
processes like the
transport and diffusion of pollutants, representing the temperature of a fluid
moving along a heated wall, or the concentration of electrons in models of
semiconductor devices, it is also a fundamental subproblem for models of
incompressible flow.

The unsteady elliptic problem can be handled in two different ways
such as ``Rothe Method'' and ``Method of Lines''. For the ``Rothe
Method'', the time variable is discretized firstly by certain time
differencing scheme to obtain a sequence of steady problems, and
each of these problems is then solved by some spatial discretization
method. For the ``Method of Lines'', the spatial variable is
discretized firstly to obtain a system of ordinary-differential
equations (ODEs) or differential-algebraic equations (DAEs), and
certain time differencing scheme is then applied to solve the above
differential equations.

The waveform relaxation (WR) methods are powerful solvers for numerically computing
the solution of ODEs or DAEs on sequential and parallel computers, which
was first introduced by Lelarasmee in \cite{LelRS-82} for simulating the behavior of
very large-scale electrical networks. Later, there are lots of expansions and
applications of his theory; see, e.g., \cite{SL-Camp-1, Jiang-09, FL-Lewis}.
The basic idea of this class of
iteration methods is to apply relaxation technique directly to the
corresponding differential equations, which can be regarded as a natural
extension of the classical relaxation methods for solving systems of linear
equations with iterating space changing from $\reals^n$ to the time-dependent
function or the waveform space.

In order to take advantage of the waveform relaxation methods for solving
the unsteady elliptic problem (\ref{uns-elliptic}), we follow the first step
of ``Method of Lines'' to discretize (\ref{uns-elliptic}) spatially
with spatial grid parameter $h$ to obtain
the so-called unsteady discrete elliptic problem as follows,
\begin{eqnarray}
\label{nsa-pd-lde} \mathcal{L}_h(x) = B\,\dot x(t) + A\, x(t) = f(t),\quad
x(0)=x_0,
\end{eqnarray}
with $B\in\ceals^{r\times r}$ being Hermitian and $A\in\ceals^{r\times r}$
being non-Hermitian positive definite, the solution $x(t)$ and the data $f(t)$
are complex vector-valued functions. It can be proved that the linear
differential operator $\mathcal{L}_h$ is non-self-adjoint positive definite
on Lebesgue square-integrable function space under suitable conditions.
Therefore, this operator $\mathcal{L}_h$ can be considered as the analogous
of non-Hermitian positive definite matrix. For systems of linear equations
related to non-Hermitian positive definite coefficient matrix, Bai, Golub
and Ng \cite{Bai-03} proposed a class of two-step iterative methods,
called the HSS method, which is designed in the spirit of the ADI iteration
technique \cite{Douglas-56} and by making use of the
natural splitting of non-Hermitian positive definite matrix, i.e., the HS
splitting. Similar to the HS splitting of matrix, we define the HS splitting
of the non-self-adjoint positive definite linear differential operator
$\mathcal{L}_h$, and design
the waveform relaxation method based on the HS spitting of operator
$\mathcal{L}_h$, i.e., the WR-HSS method, for solving the unsteady discrete
elliptic problem (\ref{nsa-pd-lde}).

The paper is organized as follows. It is started in Section
\ref{sec-hss} by reviewing the basic idea of the WR method and the
HSS method, then the framework of the WR-HSS method is described
specifically. In Section \ref{sec-wrhss-convergence}, the
convergence analysis of the WR-HSS method is given. In practical
aspect, the WR-HSS method must be implemented discretely, therefore,
the discrete-time WR-HSS method and the implementation details are
stated in Section \ref{sec-implementation-details}. The numerical
results are listed in Section \ref{sec-num-example} to show the
effectiveness of the WR-HSS method and the correctness of the
theoretical results. To end this paper, we give some concluding
remarks in Section \ref{conclusions}.

{\bf Notations:}
In order to make the meaning of $2$-norms in different spaces used
in this paper more clear, we use different notations for $2$-norms in
different spaces. Specifically, we denote $\bbb{L}_2^r(\reals)$ as
the Hilbert space consisting of
complex vector-valued functions with the inner product
\begin{eqnarray*}
(f(t),g(t))=\int_{-\infty}^{+\infty} g^*(t)f(t)\mbox{d}t
=\int_{-\infty}^{+\infty} \sum_{i=1}^r \overline{g_i(t)}f_i(t)\mbox{d}t,\
\forall\, f(t),g(t)\in \bbb{L}_2^r(\reals),
\end{eqnarray*}
where the integral is in the Lebesgue sense, and the corresponding $2$-norm is denoted
as $\|f(t)\|_{\bbb{L}}=\sqrt{(f(t),f(t))}$, $\forall\ f(t)\in \bbb{L}_2^r(\reals)$.
For convenience, we also denote
\begin{eqnarray*}
(u,v)=v^*u=\sum_i^r \overline{v_i}\,u_i,\ \forall\ u,v\in\ceals^r
\end{eqnarray*}
as the inner product of the $r$-dimensional complex vector space $\ceals^r$, and
the corresponding $2$-norm is denoted as $\|u\|_{\ceals}=\sqrt{(u,u)}$, $\forall\ u\in\ceals^r$.

\section{The WR method, HSS method and WR-HSS method}
\label{sec-hss}

In this section, we review the WR method for solving the system of linear differential
equations and the HSS method for solving the system of linear equations,
and present the WR-HSS method for solving the unsteady discrete elliptic problem
(\ref{nsa-pd-lde}).

\subsection{The WR method}
\label{sub-sec-lde-wr}

The WR method is a powerful solver for solving the system
of linear differential equations of the form (\ref{nsa-pd-lde}), i.e.,
\begin{eqnarray*}
    \mathcal{L}_h\, x(t) = f(t),
\end{eqnarray*}
which arise in abroad range of applications
in scientific/engineering computing.

By denoting
\begin{eqnarray*}
\mathcal{M}=M_B\frac{\mbox{d}}{\mbox{d}t}+M_A
\quad\mbox{and}\quad
\mathcal{N}=N_B\frac{\mbox{d}}{\mbox{d}t}+N_A
\end{eqnarray*}
with matrix splittings
\begin{eqnarray*}
B=M_B-N_B \quad\mbox{and}\quad A=M_A-N_A,
\end{eqnarray*}
we have the following operator splitting,
\begin{eqnarray*}
\mathcal{L}_h=\mathcal{M}-\mathcal{N}.
\end{eqnarray*}
The WR method is defined in the operator form
\begin{eqnarray*}
\mathcal{M}\,x^{(k+1)}(t)=\mathcal{N}\,x^{(k)}(t)+f(t),
\end{eqnarray*}
or formally written into the following fixed-point iteration form,
\begin{eqnarray*}
\,x^{(k+1)}(t)=\mathcal{K}\,x^{(k)}(t)+c(t),
\end{eqnarray*}
where $\mathcal{K}=\mathcal{M}^{-1}\mathcal{N}$ and $c=\mathcal{M}^{-1}\,f$.

The convergence theory of the WR method for the system of linear ordinary
differential equations (ODEs), i.e., the coefficient matrix $B$ being nonsingular,
has been perfectly figured out; see
\cite{wr-lcc-4,Jan-97,wr-lcc-1,wr-lcc-3,wr-lcc-5, WangBai05}.
The convergence rate of the above WR method is
\begin{eqnarray}
\label{lde-wr-convergencerate}
\rho(\mathcal{K}) = \sup_{\omega\in\reals}
\rho(\tilde{\mathcal{M}}^{-1}\tilde{\mathcal{N}}),
\end{eqnarray}
where
\begin{eqnarray*}
\tilde{\mathcal{M}} = \imath\omega M_B + M_A \quad\mbox{and}\quad
\tilde{\mathcal{N}} = \imath\omega N_B + N_A
\end{eqnarray*}
are the frequency counterpart of the operators $\mathcal{M}$ and $\mathcal{N}$.

\subsection{The HSS method}
\label{sub-sec-lae-hss}
Many applications in scientific computing lead to the following large sparse system
of linear equations
\begin{eqnarray}
\label{lae}
A\,x=b,
\end{eqnarray}
where $A\in\ceals^{r\times r}$ is non-Hermitian positive definite, and $b\in\ceals^r$.
There is a natural Hermitian/skew-Hermitian splitting (HSS) of the coefficient matrix $A$,
i.e.,
\begin{eqnarray}
\label{lae-hs-splitting}
A=H+S,
\end{eqnarray}
with
\begin{eqnarray*}
H=\frac{1}{2}(A+A^*)\quad\mbox{and}\quad S=\frac{1}{2}(A-A^*).
\end{eqnarray*}
based on the above HS splitting and motivated by the ADI iteration technique \cite{Douglas-56}, Bai, Golub and Ng proposed a class
of two-step iterative methods called the Hermitian/skew-Hermitian
splitting method; see \cite{Bai-03}.

{\bf The HSS method.} {\it Given an initial guess $x^{(0)}\in\ceals^r$,
for $k=0,1,2,\ldots$, until $\{x^{(k)}\}\subset\ceals^r$ converges, compute
\begin{eqnarray}
\label{lae-hss-iteration}
\left\{
\begin{array}{lll}
  (\alpha I+H)\,x^{(k+\frac{1}{2})} &=& (\alpha I-S)\,x^{(k)}+b,\\
  (\alpha I+S)\,x^{(k+1)} &=& (\alpha I-H)\,x^{(k+\frac{1}{2})}+b,
\end{array}
\right.
\end{eqnarray}
where $\alpha$ is a given positive constant.}

The above HSS method can be equivalently rewritten into the following
matrix-vector form
\begin{eqnarray*}
x^{(k+1)} = F(\alpha)^{-1}G(\alpha)\,x^{(k)} + F(\alpha)^{-1}\,b,\quad k=0,1,2,\ldots,
\end{eqnarray*}
where the iteration matrix results from the spitting
\begin{eqnarray*}
A=F(\alpha)-G(\alpha)
\end{eqnarray*}
of the coefficient matrix $A$ with
\begin{eqnarray*}
\left\{
\begin{array}{lll}
  F(\alpha) &=& \frac{1}{2\alpha}(\alpha I+H)(\alpha I+S),\\
  G(\alpha) &=& \frac{1}{2\alpha}(\alpha I-H)(\alpha I-S).
\end{array}
\right.
\end{eqnarray*}
The convergence property of the HSS method is described in the following theorem.
\begin{theorem}
\label{th-lae-hss-iteration}
\cite{Bai-03} Let $A\in\ceals^{r\times r}$ be a non-Hermitian positive definite matrix
with the HS splitting (\ref{lae-hs-splitting}), and $\alpha$ be a positive constant.
Then the spectral radius $\rho(F(\alpha)^{-1}G(\alpha))$ of the iteration matrix of the
HSS method is bounded by
\begin{eqnarray*}
\sigma(\alpha) = \max_{\lambda_j\in\lambda(H)} \left|\frac{\alpha-\lambda_j}
{\alpha+\lambda_j}\right|,
\end{eqnarray*}
where $\lambda(H)$ is the spectral set of the matrix $H$. Therefore, it follows that
\begin{eqnarray*}
\rho(F(\alpha)^{-1}G(\alpha))\le\sigma(\alpha)<1,\quad \forall\,\alpha>0,
\end{eqnarray*}
i.e., the HSS method is convergent.

Moreover, if $\gamma_{\min}$ and $\gamma_{\max}$ are the lower and the upper bounds
of the eigenvalues of the matrix $H$, respectively, then
\begin{eqnarray*}
\alpha^* = \arg\min_{\alpha}\left\{
\max_{\gamma_{\min}\le\lambda\le\gamma_{\max}}\left|
\frac{\alpha-\lambda_j}{\alpha+\lambda_j}
\right|
\right\}=\sqrt{\gamma_{\min}\gamma_{\max}}
\end{eqnarray*}
and
\begin{eqnarray*}
\sigma(\alpha^*) = \frac{\sqrt{\gamma_{\max}}-\gamma_{\min}}
{\sqrt{\gamma_{\max}}+\gamma_{\min}}=
\frac{\sqrt{\kappa(H)}-1}{\sqrt{\kappa(H)}+1},
\end{eqnarray*}
where $\kappa(H)$ is the spectral condition number of $H$.
\end{theorem}

The above theorem demonstrates that the HSS method is
unconditionally convergent to the unique solution of the
non-Hermitian positive definite system of linear equations
(\ref{lae}), with the same convergence rate as that of the conjugate
gradient method when it is applied to a system of linear equations
with Hermitian positive definite coefficient matrix. In addition,
the upper bound of its asymptotic convergence rate is only dependent on the
spectrum of the Hermitian part $H$, but is independent of the
spectrum of the skew-Hermitian part $S$. To learn more about the HSS
method and its variants, one can refer to references
\cite{BaiGolubLi07MC,BaiGolubNg07NLAA,BaiNg08,BaiYin05} for system
of real linear equations, references \cite{BaiBenziChen10C,BaiBenziChen11NA}
for system of complex linear equations, and references
\cite{Bai09NLAA,BaiBenziChenWang13IMA,
BaiGolub07IMA,BaiLi06,BaiPan04} for system of linear equations with
block two-by-two coefficient matrix.

\subsection{The WR-HSS method}
\label{sub-sec-wrhss}

We first consider the generalization of the HSS method to the linear
operator equation on Hilbert space. Let $\mathcal{L}$ be a linear operator defined
on Hilbert space $\bbb{V}$, and the following equation is satisfied
\begin{eqnarray}
\label{loe}
\mathcal{L}\,x = f,
\end{eqnarray}
where $f\in\bbb{V}$ is given, and $x\in\bbb{V}$ is the unknown. We denote $\mathcal{L}^*$
as the adjoint operator of $\mathcal{L}$, i.e.,
\begin{eqnarray*}
(\mathcal{L}\,u,v) = (u,\mathcal{L}^*\,v),\quad u,v\in\bbb{V},
\end{eqnarray*}
here $(\cdot,\cdot)$ is the inner product in Hilbert space $\bbb{V}$. Then we can
define the HS splitting of the linear operator $\mathcal{L}$ as
\begin{eqnarray}
\label{loe-hs-splitting}
\mathcal{L}=\mathcal{H}+\mathcal{S},
\end{eqnarray}
with
\begin{eqnarray*}
\mathcal{H}=\frac{1}{2}(\mathcal{L}+\mathcal{L}^*)\quad\mbox{and}
\quad \mathcal{S}=\frac{1}{2}(\mathcal{L}-\mathcal{L}^*).
\end{eqnarray*}
Here, we call the operators $\mathcal{H}$ and $\mathcal{S}$ as the Hermitian part
and the skew-Hermitian part of the operator $\mathcal{L}$.
Based on the above splitting, the HSS method is straightforwardly
generalized as follows.

{\bf The operatorized HSS method.} {\it Given an initial guess $x^{(0)}\in\bbb{V}$,
for $k=0,1,2,\ldots$, until $\{x^{(k)}\}\subset\bbb{V}$ converges, compute
\begin{eqnarray}
\label{loe-hss-iteration}
\left\{
\begin{array}{lll}
  (\alpha I+\mathcal{H})\,x^{(k+\frac{1}{2})} &=& (\alpha I-\mathcal{S})\,x^{(k)}+f,\\
  (\alpha I+\mathcal{S})\,x^{(k+1)} &=& (\alpha I-\mathcal{H})\,x^{(k+\frac{1}{2})}+f,
\end{array}
\right.
\end{eqnarray}
where $\alpha$ is a given positive constant.}

In the sequel, we discuss the application of the operatorized HSS mehtod (\ref{loe-hss-iteration})
to the unsteady discrete elliptic problem (\ref{nsa-pd-lde}). We
consider the solution of the unsteady discrete elliptic problem (\ref{nsa-pd-lde})
in the complex vector-valued function space $\bbb{L}_2^r(\reals)$. To do this, we need to
prolong the solution $x(t)$ and the data $f(t)$ to the whole real axis
$\reals$, and keep the notations of the prolonged functions unchanged,
i.e.,
\begin{eqnarray*}
x(t)=\left\{
\begin{array}{ll}
  x(t) & t\in\reals _+ \\
  0      & \mbox{otherwise}
\end{array}
\right.\, \mbox{and}\quad f(t)=\left\{
\begin{array}{ll}
  f(t) & t\in\reals _+ \\
  0      & \mbox{otherwise}
\end{array}
\right..
\end{eqnarray*}
Due to the definition of the Hilbert space $\bbb{L}_2^r(\reals)$, the value of a function
on a single point is no longer essential. Hence, the initial condition of the unsteady
discrete elliptic problem (\ref{nsa-pd-lde}) is ignored. Now, we can rewrite the unsteady
discrete elliptic problem (\ref{nsa-pd-lde}) in the following operator form,
\begin{eqnarray*}
\mathcal{L}_h\,x(t) = f(t), \quad x(t),f(t)\in\bbb{L}_2^r(\reals),
\end{eqnarray*}
where
\begin{eqnarray*}
\mathcal{L}_h = B\,\frac{\mbox{d}}{\mbox{d}t} + A.
\end{eqnarray*}
The adjoint of the above operator is of the form
\begin{eqnarray*}
\mathcal{L}^*_h = -B^*\,\frac{\mbox{d}}{\mbox{d}t} + A^*,
\end{eqnarray*}
which can be verified as
\begin{eqnarray*}
(\mathcal{L}_h\,u(t),v(t))
& = & \int_{-\infty}^{+\infty} v^*(t)(B\,\dot u(t) + A\,u(t))\mbox{d}t\\
& = & \int_{-\infty}^{+\infty} v^*(t)B\, \mbox{d}u(t) +
      \int_{-\infty}^{+\infty} v^*(t) A\,u(t) \mbox{d}t\\
& = & v^*(t)B\,u(t)|_{-\infty}^{+\infty} -
      \int_{-\infty}^{+\infty} (\mbox{d}v^*(t))B\,u(t) + (u(t),A^*\,v(t))\\
& = & (u(t),-B^*\,\dot v(t)) + (u(t),A^*\,v(t))\\
& = & (u(t),\mathcal{L}^*_h\,v(t)),\quad \forall\,u(t),v(t)\in\bbb{L}_2^r(\reals),
\end{eqnarray*}
the derivatives $\dot u(t)$ and $\dot v(t)$ are taken in the sense of distribution.
Then we have the HS splitting of the operator $\mathcal{L}_h$ as follows
\begin{eqnarray*}
\mathcal{L}_h = \mathcal{H}_h + \mathcal{S}_h,
\end{eqnarray*}
where
\begin{eqnarray*}
\mathcal{H}_h = \frac{1}{2}(\mathcal{L}_h+\mathcal{L}^*_h) =
\frac{1}{2}(B-B^*)\frac{\mbox{d}}{\mbox{d}t} + H
\end{eqnarray*}
and
\begin{eqnarray*}
\mathcal{S}_h = \frac{1}{2}(\mathcal{L}_h-\mathcal{L}^*_h) =
\frac{1}{2}(B+B^*)\frac{\mbox{d}}{\mbox{d}t} + S,
\end{eqnarray*}
here $H$ and $S$ are the Hermitian part and the skew-Hermitian part of
the coefficient matrix $A$ respectively. Since the matrix $B$ is Hermitian,
the above expression can be simplified as
\begin{eqnarray*}
\mathcal{H}_h = H \quad\mbox{and}\quad \mathcal{S}_h =
B\frac{\mbox{d}}{\mbox{d}t} + S.
\end{eqnarray*}
Based on the HS splitting of the operator $\mathcal{L}_h$, the operatorized
HSS method (\ref{loe-hss-iteration}) becomes the following iteration
scheme.

{\bf The WR-HSS method.} {\it Given an initial guess $x^{(0)}\in\bbb{L}_2^r(\reals)$,
for $k=0,1,2,\ldots$, until $\{x^{(k)}\}\subset\bbb{L}_2^r(\reals)$ converges, compute
\begin{eqnarray}
\label{lde-wrhss-iteration-1}
\left\{
\begin{array}{lll}
  (\alpha I+\mathcal{H}_h)\,x^{(k+\frac{1}{2})} &=& (\alpha I-\mathcal{S}_h)\,x^{(k)}+f,\\
  (\alpha I+\mathcal{S}_h)\,x^{(k+1)} &=& (\alpha I-\mathcal{H}_h)\,x^{(k+\frac{1}{2})}+f,
\end{array}
\right.
\end{eqnarray}
or equivalently,
\begin{eqnarray}
\label{lde-wrhss-iteration-2}
\left\{
\begin{array}{lll}
  (\alpha I+H)\,x^{(k+\frac{1}{2})} &=& (\alpha I-B\,\frac{\mbox{d}}{\mbox{d}t}-S)\,x^{(k)}+f,\\
  (\alpha I+B\,\frac{\mbox{d}}{\mbox{d}t}+S)\,x^{(k+1)} &=& (\alpha I-H)\,x^{(k+\frac{1}{2})}+f,
\end{array}
\right.
\end{eqnarray}
where $\alpha$ is a given positive constant.}

The above WR-HSS method can be formally rewritten into the following
fixed-point iteration form
\begin{eqnarray*}
x^{(k+1)} = \mathcal{K}_{\mbox{\tiny WR-HSS}}\,x^{(k)} + c,\quad k=0,1,2,\ldots,
\end{eqnarray*}
with $\mathcal{K}_{\mbox{\tiny WR-HSS}}=\mathcal{F}(\alpha)^{-1}\mathcal{G}(\alpha)$ and
$c=\mathcal{F}(\alpha)^{-1}f$.
The iteration operator of the above iteration scheme results from the splitting
\begin{eqnarray*}
\mathcal{L}_h=\mathcal{F}(\alpha)-\mathcal{G}(\alpha)
\end{eqnarray*}
with
\begin{eqnarray*}
\left\{
\begin{array}{lll}
  \mathcal{F}(\alpha) &=& \frac{1}{2\alpha}(\alpha I+\mathcal{H}_h)(\alpha I+\mathcal{S}_h),\\
  \mathcal{G}(\alpha) &=& \frac{1}{2\alpha}(\alpha I-\mathcal{H}_h)(\alpha I-\mathcal{S}_h).
\end{array}
\right.
\end{eqnarray*}

We end this section with two remarks.
\begin{remark}
\label{rema-ohss}
$\,$
\begin{itemize}
    \item If we reverse the roles of the operators $\mathcal{H}_h$ and $\mathcal{S}_h$ in the
    WR-HSS method (\ref{lde-wrhss-iteration-1}), then we obtain the following iteration scheme,
    \begin{eqnarray} \label{lde-wrhss-iteration-var-1}
        \left\{
        \begin{array}{lll}
        (\alpha I+\mathcal{S}_h)\,x^{(k+\frac{1}{2})} &=& (\alpha I-\mathcal{H}_h)\,x^{(k)}+f,\\
        (\alpha I+\mathcal{H}_h)\,x^{(k+1)} &=& (\alpha I-\mathcal{S}_h)\,x^{(k+\frac{1}{2})}+f,
        \end{array}
        \right.
    \end{eqnarray}
    or equivalently,
    \begin{eqnarray*}
        \left\{
        \begin{array}{lll}
        (\alpha I+B\,\frac{\mbox{d}}{\mbox{d}t}+S)\,x^{(k+\frac{1}{2})} &=& (\alpha I-H)\,x^{(k)}+f,\\
        (\alpha I+H)\,x^{(k+1)} &=& (\alpha I-B\,\frac{\mbox{d}}{\mbox{d}t}-S)\,x^{(k+\frac{1}{2})}+f,
        \end{array}
        \right.
    \end{eqnarray*}
    the convergence properties of the above iteration scheme are similar to those of the WR-HSS method.
    \item If we only reverse the roles of the matrices $H$ and $S$ in the WR-HSS method
    (\ref{lde-wrhss-iteration-2}), then we have the following iteration scheme
    \begin{eqnarray}
    \label{lde-wrhss-bad}
        \left\{
        \begin{array}{lll}
        (\alpha I+S)\,x^{(k+\frac{1}{2})} &=& (\alpha I-B\,\frac{\mbox{d}}{\mbox{d}t}-H)\,x^{(k)}+f,\\
        (\alpha I+B\,\frac{\mbox{d}}{\mbox{d}t}+H)\,x^{(k+1)} &=& (\alpha I-S)\,x^{(k+\frac{1}{2})}+f,
        \end{array}
        \right.
    \end{eqnarray}
    this is just a slight change of the original WR-HSS method, but the convergence rate
    of the so modified WR-HSS method is much slower. The reason is that the
    operator splitting related to the above iteration scheme is no longer a HS splitting.
\end{itemize}
\end{remark}

\section{Convergence analysis of the WR-HSS method}
\label{sec-wrhss-convergence}

In this section, we study the convergence properties of the WR-HSS method.
The convergence analysis is carried out with the help of the Fourier transform.
Firstly, we choose the definition of the Fourier transform of a function $v(t)$
as follows
\begin{eqnarray*}
\mathcal{F}v(t)=\frac{1}{\sqrt{2\pi}}\int_{-\infty}^{+\infty}e^{-\imath\omega t}v(t)\mbox{d}t
=\tilde v(\omega),
\end{eqnarray*}
then the above Fourier transform is a unitary operator mapping the space $\bbb{L}_2^r(\reals)$
into itself; see \cite{Hochstadt-89}. In addition, we assume the existence of the Fourier
transforms of the solution $x(t)$, the data $f(t)$, the iterate $x^{(k)}$ and the intermediate
iterate $x^{(k+\frac{1}{2})}$.

We transform the unsteady discrete elliptic problem in time domain (\ref{nsa-pd-lde}) to
its counterpart in frequency domain as follows
\begin{eqnarray}
\label{nsa-pd-lde-frequency}
\tilde{\mathcal{L}}_h\, \tilde x(\omega) = \tilde f(\omega),
\end{eqnarray}
where $\tilde{\mathcal{L}}_h=\imath \omega B + A$, $\tilde x(\omega) = \mathcal{F}x(t)$
and $\tilde f(\omega) = \mathcal{F}f(t)$. Obviously, the HS splitting of the operator
$\tilde{\mathcal{L}}_h$ is given by
\begin{eqnarray*}
\tilde{\mathcal{L}}_h = \tilde{\mathcal{H}}_h + \tilde{\mathcal{S}}_h
\end{eqnarray*}
with $\tilde{\mathcal{H}}_h = H$ and $\tilde{\mathcal{S}}_h = \imath\omega B + S$.
Due to the properties of matrices $B$ and $A$, we remark that $\tilde{\mathcal{L}}_h$
is non-Hermitian positive definite for any given frequency $\omega\in\reals$.

Similarly, the application of the Fourier transform to the WR-HSS method
(\ref{lde-wrhss-iteration-1}) leads to its counterpart iteration scheme in
frequency domain, i.e.,
\begin{eqnarray}
\label{lde-wrhss-frequency}
\left\{
\begin{array}{lll}
  (\alpha I+\tilde{\mathcal{H}}_h)\,\tilde x^{(k+\frac{1}{2})} &=& (\alpha I-\tilde{\mathcal{S}}_h)\,\tilde x^{(k)}+\tilde f,\\
  (\alpha I+\tilde{\mathcal{S}}_h)\,\tilde x^{(k+1)} &=& (\alpha I-\tilde{\mathcal{H}}_h)\,\tilde x^{(k+\frac{1}{2})}+\tilde f.
\end{array}
\right.
\end{eqnarray}
For any given frequency $\omega\in\reals$, the above iteration scheme is just the HSS
method of the system of linear equations (\ref{nsa-pd-lde-frequency}), which
can be rewritten into the following
matrix-vector form
\begin{eqnarray*}
\tilde x^{(k+1)} = \tilde K_{\mbox{\tiny WR-HSS}}\,\tilde x^{(k)} + \tilde c
\end{eqnarray*}
with $\tilde K_{\mbox{\tiny WR-HSS}} = \tilde F(\alpha)^{-1}\tilde G(\alpha)$
and $\tilde c = \tilde F(\alpha)^{-1}\,\tilde f$.
The iteration matrix of the above iteration scheme results from the splitting
\begin{eqnarray*}
\tilde{\mathcal{L}}_h=\tilde F(\alpha)-\tilde G(\alpha)
\end{eqnarray*}
of the matrix $\tilde{\mathcal{L}}_h$ with
\begin{eqnarray*}
\left\{
\begin{array}{lll}
  \tilde F(\alpha) &=& \frac{1}{2\alpha}(\alpha I+\tilde{\mathcal{H}}_h)(\alpha I+\tilde{\mathcal{S}}_h),\\
  \tilde G(\alpha) &=& \frac{1}{2\alpha}(\alpha I-\tilde{\mathcal{H}}_h)(\alpha I-\tilde{\mathcal{S}}_h),
\end{array}
\right.
\end{eqnarray*}
where $\tilde F(\alpha)$ and $\tilde G(\alpha)$ are just the frequency counterparts of the operators
$\mathcal{F}(\alpha)$ and $\mathcal{G}(\alpha)$ respectively.

{\bf Convergence in frequency domain.}
According to the theory of the HSS method (i.e., Theorem \ref{th-lae-hss-iteration}),
we have the convergence property of
the WR-HSS method in frequency domain (\ref{lde-wrhss-frequency}).
\begin{theorem}
\label{th-lde-wrhss-frequency}
Consider the WR-HSS method
in frequency domain (\ref{lde-wrhss-frequency}) for the unsteady discrete
elliptic problem in frequency domain (\ref{nsa-pd-lde-frequency}).
Let $\alpha$ be a positive constant.
Then the spectral radius $\rho(\tilde K_{\mbox{\tiny WR-HSS}})$ of
the iteration matrix is bounded by
\begin{eqnarray*}
\sigma(\alpha) = \max_{\lambda_j\in\lambda(\tilde{\mathcal{H}}_h)} \left|\frac{\alpha-\lambda_j}
{\alpha+\lambda_j}\right|,
\end{eqnarray*}
where $\lambda(\tilde{\mathcal{H}}_h)$ is the spectral set of the matrix
$\tilde{\mathcal{H}}_h = H$. Therefore, it follows that
\begin{eqnarray*}
\rho(\tilde K_{\mbox{\tiny WR-HSS}})\le\sigma(\alpha)<1,\quad \forall\,\alpha>0,
\end{eqnarray*}
i.e., the WR-HSS method in frequency domain (\ref{lde-wrhss-frequency})
is convergent for any given frequency $\omega\in\reals$.

Moreover, if $\gamma_{\min}$ and $\gamma_{\max}$ are the lower and the
upper bounds of the eigenvalues of the matrix $\tilde{\mathcal{H}}_h = H$,
respectively, then
\begin{eqnarray*}
\alpha^* = \arg\min_{\alpha}\left\{
\max_{\gamma_{\min}\le\lambda\le\gamma_{\max}}\left|
\frac{\alpha-\lambda_j}{\alpha+\lambda_j}
\right|
\right\}=\sqrt{\gamma_{\min}\gamma_{\max}}
\end{eqnarray*}
and
\begin{eqnarray*}
\sigma(\alpha^*) = \frac{\sqrt{\gamma_{\max}}-\gamma_{\min}}
{\sqrt{\gamma_{\max}}+\gamma_{\min}}=
\frac{\sqrt{\kappa(\tilde{\mathcal{H}}_h)}-1}{\sqrt{\kappa(\tilde{\mathcal{H}}_h)}+1},
\end{eqnarray*}
where $\tilde{\mathcal{H}}_h$ is the spectral condition number of
$\tilde{\mathcal{H}}_h$.
\end{theorem}

\begin{remark}
\label{rema-wrhss-convergence} $\,$
\begin{itemize}
    \item Theorem \ref{th-lde-wrhss-frequency} shows that the asymptotic convergence
    rate of the WR-HSS method in frequency domain (\ref{lde-wrhss-frequency}) is bounded
    by the positive real function $\sigma(\alpha)$, which is also an upper bound of the
    asymptotic convergence rate of the WR-HSS method in time domain (\ref{lde-wrhss-iteration-1})
    according to the result (\ref{lde-wr-convergencerate}), i.e.,
    \begin{eqnarray*}
    \rho(\mathcal{K}_{\mbox{\tiny WR-HSS}}) &=& \sup_{\omega\in\reals}
    \rho(\tilde F(\alpha)^{-1}\tilde G(\alpha)) \\
    &=& \sup_{\omega\in\reals}
    \rho((\alpha I+\tilde{\mathcal{S}}_h)^{-1}(\alpha I+\tilde{\mathcal{H}}_h)^{-1}
    (\alpha I-\tilde{\mathcal{H}}_h)(\alpha I-\tilde{\mathcal{S}}_h))\\
    &=& \sup_{\omega\in\reals}
    \rho((\alpha I+\tilde{\mathcal{H}}_h)^{-1}(\alpha I-\tilde{\mathcal{H}}_h)
    (\alpha I-\tilde{\mathcal{S}}_h)(\alpha I+\tilde{\mathcal{S}}_h)^{-1})\\
    &\le& \sup_{\omega\in\reals}
    \|(\alpha I+\tilde{\mathcal{H}}_h)^{-1}(\alpha I-\tilde{\mathcal{H}}_h)
    (\alpha I-\tilde{\mathcal{S}}_h)(\alpha I+\tilde{\mathcal{S}}_h)^{-1}\|_2\\
    &\le& \sup_{\omega\in\reals}
    \|(\alpha I+\tilde{\mathcal{H}}_h)^{-1}(\alpha I-\tilde{\mathcal{H}}_h)\|_2
    \|(\alpha I-\tilde{\mathcal{S}}_h)(\alpha I+\tilde{\mathcal{S}}_h)^{-1}\|_2\\
    &=& \|(\alpha I+H)^{-1}(\alpha I-H)\|_2 \\
    &=& \sigma(\alpha).
    \end{eqnarray*}
    \item According to the analysis in \cite{BertGolub05}, for a steady convection
    dominated elliptic problem, the positive real function $\sigma(\alpha)$ is close to one, but the
    asymptotic convergence rate
    of the corresponding HSS method is far less than the positive real function
    $\sigma(\alpha)$. Then, for an unsteady convection
    dominated elliptic problem, the positive real function
    $\sigma(\alpha)$ is also close to one and the asymptotic convergence
    rate $\rho(\mathcal{K}_{\mbox{\tiny WR-HSS}})$ of the WR-HSS method
    might satisfy the following inequality
    \begin{eqnarray*}
    \rho(\mathcal{K}_{\mbox{\tiny WR-HSS}})\ll\sigma(\alpha)\approx 1,
    \end{eqnarray*}
    which means that
    the WR-HSS method could perform much better than the upper bound $\sigma(\alpha)$
    can reveal in convection dominated cases.
\end{itemize}
\end{remark}

Now, we give an explanation of why the convergence property of the iteration
scheme (\ref{lde-wrhss-bad}) in Remark \ref{rema-ohss} is not as good as
the WR-HSS method and its variant (\ref{lde-wrhss-iteration-var-1}).
We apply the Fourier transform to the iteration scheme (\ref{lde-wrhss-bad})
and obtain its counterpart in frequency domain as follows
\begin{eqnarray} \label{lde-wrhss-bad-frequency}
\left\{
\begin{array}{lll}
  (\alpha I+S)\,\tilde x^{(k+\frac{1}{2})} &=& (\alpha I-\imath\omega B - H)\,\tilde x^{(k)}+\tilde f,\\
  (\alpha I+\imath\omega B + H)\,\tilde x^{(k+1)} &=& (\alpha I-S)\,\tilde x^{(k+\frac{1}{2})}+\tilde f.
\end{array}
\right.
\end{eqnarray}
Obviously, the matrix $S$ is skew-Hermitian, but the matrix $\imath\omega B + H$
is not Hermitian, which means that the above two matrices do not compose a HS
splitting of the operator $\tilde{\mathcal{L}}_h$ of the unsteady discrete
elliptic problem (\ref{nsa-pd-lde-frequency}) in frequency domain. Therefore,
the iteration scheme (\ref{lde-wrhss-bad-frequency}) does not have the convergence
property of the frequency counterpart (\ref{lde-wrhss-frequency}) of the WR-HSS method,
or say, the iteration scheme (\ref{lde-wrhss-bad}) does not have the convergence
property of the WR-HSS method. More specifically, direct computation leads to
the iteration matrix of the iteration scheme (\ref{lde-wrhss-bad-frequency}),
i.e.,
\begin{eqnarray*}
\mathcal{G}(\alpha;\omega)=
(\alpha I+\imath\omega B + H)^{-1}(\alpha I-S)(\alpha I+S)^{-1}(\alpha I-\imath\omega B - H),
\end{eqnarray*}
which is similar to the following matrix
\begin{eqnarray*}
\hat{\mathcal{G}}(\alpha;\omega)=
(\alpha I-\imath\omega B - H)(\alpha I+\imath\omega B + H)^{-1}(\alpha I-S)(\alpha I+S)^{-1}.
\end{eqnarray*}
Here matrix $\alpha I+S$ and matrix $\alpha I+\imath\omega B + H$ are nonsingular for
any positive constant $\alpha$. Then we have
\begin{eqnarray*}
\rho(\mathcal{G}(\alpha;\omega))
    & = &\rho(\hat{\mathcal{G}}(\alpha;\omega))\\
    & \le & \|(\alpha I-\imath\omega B - H)(\alpha I+\imath\omega B + H)^{-1}(\alpha I-S)(\alpha I+S)^{-1}\|_2\\
    & \le & \|(\alpha I-\imath\omega B - H)(\alpha I+\imath\omega B + H)^{-1}\|_2\|(\alpha I-S)(\alpha I+S)^{-1}\|_2.
\end{eqnarray*}
Since $Q(\alpha)=(\alpha I-S)(\alpha I+S)^{-1}$ is the Cayley transform of the skew-Hermitian matrix $S$,
it means that $Q(\alpha)$ is a unitary matrix. Therefore, $\|(\alpha I-S)(\alpha I+S)^{-1}\|_2=1$.
If we assume that the matrices $B$ and $H$ are commutative (the fact is just the case or can be
equivalently transformed to the case of such kind in most of the time), then the matrix
$\imath\omega B + H$ is normal. It follows that
\begin{eqnarray*}
\rho(\mathcal{G}(\alpha;\omega)) \le
\|(\alpha I-\imath\omega B - H)(\alpha I+\imath\omega B + H)^{-1}\|_2 =
\hat{\sigma}(\alpha;\omega),
\end{eqnarray*}
with
\begin{eqnarray*}
\hat{\sigma}(\alpha;\omega)=\max_{\lambda_j(\omega)\in\lambda(\imath\omega B+H)}
\left|\frac{\alpha-\lambda_j(\omega)}{\alpha+\lambda_j(\omega)}\right|.
\end{eqnarray*}
Since matrices $B$ and $H$ are Hermitian positive definite, the real part of each
$\lambda_j(\omega)$ satisfies $0<\Re(\lambda_j(\omega))<+\infty$, and the imaginary
part of each $\lambda_j(\omega)$ satisfies
$\lim _{\omega\rightarrow\infty} \Im(\lambda_j(\omega))=\infty$. Hence, we have
\begin{eqnarray*}
\rho(\mathcal{G}(\alpha;\omega)) \le
\hat{\sigma}(\alpha;\omega)<1 \quad\mbox{and}\quad
\sup _{\omega\in\reals} \hat{\sigma}(\alpha;\omega) = 1.
\end{eqnarray*}
These demonstrate that, the convergence rate of the iteration scheme
(\ref{lde-wrhss-bad-frequency}) is less than one which guarantee the
convergence of the iteration schemes (\ref{lde-wrhss-bad-frequency})
and (\ref{lde-wrhss-bad}), but the supremum of the upper bound
$\hat{\sigma}(\alpha;\omega)$ of the convergence rate with respect
to the frequency $\omega$ is equal to one which means that the convergence
might be very slow.

{\bf Convergence in time domain.}
Since the WR-HSS method is an iterative method in time domain rather than
in frequency domain, the convergence analysis in frequency domain does not give a full
picture of the convergence behavior of the WR-HSS method. In addition,
Remark \ref{rema-wrhss-convergence}
proved that the factor $\sigma(\alpha)$ is an upper bound of the asymptotic
convergence rate of the WR-HSS method, but it gives no implication on the
contraction property of the WR-HSS method in each iteration.
Therefore, it is necessary to discuss the contraction property of the WR-HSS method
in time domain. In fact, we can prove that the factor $\sigma(\alpha)$ is also
an upper bound of the contraction factor of the WR-HSS method on each iteration.

We introduce a norm of vector-valued function in time domain as
$\interleave\cdot\interleave _t = \|(\alpha I + \mathcal{S}_h)\cdot\|_{\bbb{L}}$,
and denote $\bbb{V}_t$ as the completion of the linear span of the set
$\{v(t) \, | \, v(t)\in\bbb{L}_2^r(\reals)\,\mbox{and}\,\interleave v(t) \interleave _t < +\infty \}$.
Then $\bbb{V}_t$ is a Banach space under the norm $\interleave\cdot\interleave _t$.
In addition, we also introduce a norm of vector-valued function in frequency domain as
$\interleave\cdot\interleave _{\omega} = \|(\alpha I + \tilde{\mathcal{S}}_h)\cdot\|_{\bbb{L}}$,
and denote $\bbb{V}_{\omega}$ as the completion of the linear span of the set
$\{\tilde v(\omega) \, | \, \tilde v(\omega)\in\bbb{L}_2^r(\reals)\,\mbox{and}\,\interleave \tilde v(\omega) \interleave _{\omega} < +\infty \}$.
Then $\bbb{V}_{\omega}$ is a Banach space under the norm $\interleave\cdot\interleave _{\omega}$.

Based on the above notations and definitions, we have the following lemma.
\begin{lemma}
\label{lemma-norm}
If $\tilde v$ is the Fourier transform of $v$, then $\tilde v \in \bbb{V}_{\omega}$
if and only if $v \in \bbb{V}_{t}$.
\end{lemma}

{\em Proof.} By direct computation, we obtain the following fact
\begin{eqnarray*}
  \interleave\tilde v \interleave _{\omega} &=& \|(\alpha I + \tilde{\mathcal{S}}_h)\tilde v\|_{\bbb{L}} \\
  &=& \| (\alpha I + \imath\omega B + S)\tilde v\|_{\bbb{L}} \\
  &=& \| \mathcal{F}^{-1}\{(\alpha I + \imath\omega B + S)\tilde v\}\|_{\bbb{L}} \\
  &=& \| (\alpha I + B \frac{\mbox{d}}{\mbox{d}t} + S) v \|_{\bbb{L}} \\
  &=& \| (\alpha I + \mathcal{S}_h) v \|_{\bbb{L}} \\
  &=& \interleave v \interleave _{t}.
\end{eqnarray*}
Hence, $\tilde v \in \bbb{V}_{\omega}$
if and only if $v \in \bbb{V}_{t}$.
\hfill $\square$

Assume that $f\in\bbb{L}_2^r(\reals)$ and $x_{\star}\in\bbb{L}_2^r(\reals)$
satisfies the following linear operator equation
\begin{eqnarray*}
  \mathcal{L}_h\, x_{\star} = f,
\end{eqnarray*}
then we have
\begin{eqnarray*}
  (\alpha I + \mathcal{S}_h)\, x_{\star} = (\alpha I - \mathcal{H}_h)\, x_{\star} + f,
\end{eqnarray*}
or equivalently,
\begin{eqnarray*}
  (\alpha I + \mathcal{S}_h)\, x_{\star} = (\alpha I - H)\, x_{\star} + f
\end{eqnarray*}
for the fact that $\mathcal{H}_h=H$. Hence, we have
\begin{eqnarray*}
  \|(\alpha I + \mathcal{S}_h)\, x_{\star}\|_{\bbb{L}} &=& \|(\alpha I - H)\, x_{\star} + f\|_{\bbb{L}} \\
  &\le& \|(\alpha I - H)\, x_{\star}\|_{\bbb{L}} + \|f\|_{\bbb{L}},
\end{eqnarray*}
i.e.,
\begin{eqnarray*}
  \interleave x_{\star}\interleave _t \le \|(\alpha I - H)\, x_{\star}\|_{\bbb{L}} + \|f\|_{\bbb{L}}.
\end{eqnarray*}
Since $f,x_{\star}\in\bbb{L}_2^r(\reals)$, we have
\begin{eqnarray*}
  \interleave x_{\star}\interleave _t < +\infty,
\end{eqnarray*}
which implies that $x_{\star}\in\bbb{V}_t$. According to the Lemma \ref{lemma-norm}, the above
inequality also leads to the fact that $\tilde{x}_{\star}\in\bbb{V}_{\omega}$.

Under suitable conditions, we can prove that the WR-HSS method (\ref{lde-wrhss-iteration-1})
and its frequency domain counterpart (\ref{lde-wrhss-frequency}) are closed in Banach space $\bbb{V}_t$
and Banach space $\bbb{V}_{\omega}$ respectively. Let the initial guess of the WR-HSS method
(\ref{lde-wrhss-iteration-1}) satisfies $x^{(0)}\in\bbb{V}_t\cap\bbb{L}_2^r(\reals)$.
In addition, we assume that the $(k)$-th iterate satisfies $x^{(k)}\in\bbb{V}_t\cap\bbb{L}_2^r(\reals)$,
then we prove the $(k+1)$-th iterate also satisfies $x^{(k+1)}\in\bbb{V}_t\cap\bbb{L}_2^r(\reals)$.
There are two half steps in the $(k+1)$-th iteration of the WR-HSS method, i.e.,
\begin{eqnarray*}
\left\{
\begin{array}{lll}
  (\alpha I+\mathcal{H}_h)\,x^{(k+\frac{1}{2})} &=& (\alpha I-\mathcal{S}_h)\,x^{(k)}+f,\\
  (\alpha I+\mathcal{S}_h)\,x^{(k+1)} &=& (\alpha I-\mathcal{H}_h)\,x^{(k+\frac{1}{2})}+f,
\end{array}
\right.
\end{eqnarray*}
or equivalently,
\begin{eqnarray*}
\left\{
\begin{array}{lll}
  (\alpha I+H)\,x^{(k+\frac{1}{2})} &=& (\alpha I-\mathcal{S}_h)\,x^{(k)}+f,\\
  (\alpha I+\mathcal{S}_h)\,x^{(k+1)} &=& (\alpha I-H)\,x^{(k+\frac{1}{2})}+f,
\end{array}
\right.
\end{eqnarray*}
for the fact that $\mathcal{H}_h=H$. From the first half step of the above iteration,
we can straightforwardly obtain that
\begin{eqnarray*}
  \|(\alpha I + H)\, x^{(k+\frac{1}{2})}\|_{\bbb{L}}
  &=& \|(\alpha I-\mathcal{S}_h)\,x^{(k)}+f\|_{\bbb{L}} \\
  &=& \|(-\alpha I-\mathcal{S}_h)\,x^{(k)}+2\alpha x^{(k)}+f\|_{\bbb{L}} \\
  &\le& \interleave x^{(k)} \interleave _t + 2\alpha\|x^{(k)}\|_{\bbb{L}}+\|f\|_{\bbb{L}} \\
  &<& +\infty,
\end{eqnarray*}
where the last inequality is because the assumptions
$x^{(k)}\in\bbb{V}_t\cap\bbb{L}_2^r(\reals)$ and $f\in\bbb{L}_2^r(\reals)$.
Since the matrix $H$ is Hermitian positive definite and
the real number $\alpha$ is positive, we have
\begin{eqnarray*}
  \|x^{(k+\frac{1}{2})}\|_{\bbb{L}} < +\infty,
\end{eqnarray*}
which implies that
\begin{eqnarray*}
  \|(\alpha I - H)\, x^{(k+\frac{1}{2})}\|_{\bbb{L}}
  < +\infty.
\end{eqnarray*}
From the second half step of the $(k+1)$-th iteration of the WR-HSS method,
we obtain
\begin{eqnarray*}
  \|(\alpha I+\mathcal{S}_h)\,x^{(k+1)}\|_{\bbb{L}}
  &=& \|(\alpha I-H)\,x^{(k+\frac{1}{2})}+f\|_{\bbb{L}}\\
  &\le& \|(\alpha I-H)\,x^{(k+\frac{1}{2})}\|_{\bbb{L}}+\|f\|_{\bbb{L}} \\
  &<& +\infty,
\end{eqnarray*}
i.e.,
\begin{eqnarray*}
  \interleave x^{(k+1)}\interleave_t
  < +\infty.
\end{eqnarray*}
Since $\alpha$ is positive and $\mathcal{S}_h$ is skew-Hermitian, we can prove that
$x^{(k+1)}$ is square integrable based on the above fact.
Therefore, we have $x^{(k+1)}\in\bbb{V}_t\cap\bbb{L}_2^r(\reals)$.
According to Lemma \ref{lemma-norm},
the frequency counterpart of $x^{(k+1)}$ also satisfies
$\tilde{x}^{(k+1)}\in\bbb{V}_{\omega}\cap\bbb{L}_2^r(\reals)$.

Now, we derive the upper bound of the contraction factor of the
WR-HSS method in time domain, i.e., in the Banach space $\bbb{V}_t$.
Based on Lemma \ref{lemma-norm}, we have
\begin{eqnarray*}
  \interleave x^{(k+1)}-x_{\star}\interleave_t
  &=& \interleave \tilde{x}^{(k+1)}-\tilde{x}_{\star}\interleave_{\omega} \\
  &=& \|(\alpha I + \tilde{\mathcal{S}}_h)(\tilde{x}^{(k+1)}-\tilde{x}_{\star})\|_{\bbb{L}} \\
  &=& \|(\alpha I + \tilde{\mathcal{S}}_h)\tilde F(\alpha)^{-1}\tilde G(\alpha)(\tilde{x}^{(k)}-\tilde{x}_{\star})\|_{\bbb{L}} \\
  &=& \left(\int _{-\infty}^{+\infty} \|(\alpha I + \tilde{\mathcal{S}}_h)\tilde F(\alpha)^{-1}\tilde G(\alpha)(\tilde{x}^{(k)}-\tilde{x}_{\star})\|_{\ceals}^2\, \mbox{d}\omega \right)^{\frac{1}{2}} \\
  &=& \left(\int _{-\infty}^{+\infty} \|(\alpha I + \tilde{\mathcal{S}}_h)\tilde F(\alpha)^{-1}\tilde G(\alpha)(\alpha I + \tilde{\mathcal{S}}_h)^{-1}\right. \\
  & & \quad\quad\quad\left. (\alpha I + \tilde{\mathcal{S}}_h)(\tilde{x}^{(k)}-\tilde{x}_{\star})\|_{\ceals}^2\, \mbox{d}\omega \right)^{\frac{1}{2}} \\
  &\le& \left(\int _{-\infty}^{+\infty} \|(\alpha I + \tilde{\mathcal{S}}_h)\tilde F(\alpha)^{-1}\tilde G(\alpha)(\alpha I + \tilde{\mathcal{S}}_h)^{-1}\|_{\ceals}^2\, \right. \\
  & & \quad\quad\quad\left.\|(\alpha I + \tilde{\mathcal{S}}_h)(\tilde{x}^{(k)}-\tilde{x}_{\star})\|_{\ceals}^2\, \mbox{d}\omega \right)^{\frac{1}{2}} \\
  &=& \left(\int _{-\infty}^{+\infty} \|(\alpha I + \tilde{\mathcal{H}}_h)^{-1}(\alpha I - \tilde{\mathcal{H}}_h)(\alpha I - \tilde{\mathcal{S}}_h)(\alpha I + \tilde{\mathcal{S}}_h)^{-1}\|_{\ceals}^2\,\right. \\
  & & \quad\quad\quad\left.\|(\alpha I + \tilde{\mathcal{S}}_h)(\tilde{x}^{(k)}-\tilde{x}_{\star})\|_{\ceals}^2\, \mbox{d}\omega \right)^{\frac{1}{2}} \\
  &\le& \left(\int _{-\infty}^{+\infty} \|(\alpha I + \tilde{\mathcal{H}}_h)^{-1}(\alpha I - \tilde{\mathcal{H}}_h)\|_{\ceals}^2\, \|(\alpha I + \tilde{\mathcal{S}}_h)(\tilde{x}^{(k)}-\tilde{x}_{\star})\|_{\ceals}^2\, \mbox{d}\omega \right)^{\frac{1}{2}} \\
  &=& \sigma(\alpha)\left(\int _{-\infty}^{+\infty} \|(\alpha I + \tilde{\mathcal{S}}_h)(\tilde{x}^{(k)}-\tilde{x}_{\star})\|_{\ceals}^2\, \mbox{d}\omega \right)^{\frac{1}{2}} \\
  &=& \sigma(\alpha)\, \|(\alpha I + \tilde{\mathcal{S}}_h)(\tilde{x}^{(k)}-\tilde{x}_{\star})\|_{\bbb{L}} \\
  &=& \sigma(\alpha)\, \interleave \tilde{x}^{(k)}-\tilde{x}_{\star}\interleave_{\omega} \\
  &=& \sigma(\alpha)\, \interleave x^{(k)}-x_{\star}\interleave_t,
\end{eqnarray*}
where the factor $\sigma(\alpha)$ is of the following form
\begin{eqnarray*}
\sigma(\alpha)
  = \|(\alpha I + \tilde{\mathcal{H}}_h)^{-1}(\alpha I - \tilde{\mathcal{H}}_h)\|_{\ceals}
  = \max_{\lambda_j\in\lambda(\tilde{\mathcal{H}}_h)}
\left|\frac{\alpha-\lambda_j} {\alpha+\lambda_j}\right|.
\end{eqnarray*}
The above factor $\sigma(\alpha)$ is just an upper bound of the contraction factor of
the WR-HSS method in time domain (\ref{lde-wrhss-iteration-1}) under the norm
$\interleave \cdot \interleave_t$. Similarly to the analysis in \cite{Bai-03},
we can determine the optimal $\alpha$ to minimize the factor $\sigma(\alpha)$.
If $\gamma_{\min}$ and $\gamma_{\max}$ are the lower and
the upper bounds of the eigenvalues of the matrix
$\tilde{\mathcal{H}}_h = H$, respectively, then
\begin{eqnarray*}
\alpha^* = \arg\min_{\alpha}\left\{
\max_{\gamma_{\min}\le\lambda\le\gamma_{\max}}\left|
\frac{\alpha-\lambda_j}{\alpha+\lambda_j} \right|
\right\}=\sqrt{\gamma_{\min}\gamma_{\max}}
\end{eqnarray*}
and
\begin{eqnarray*}
\sigma(\alpha^*) = \frac{\sqrt{\gamma_{\max}}-\gamma_{\min}}
{\sqrt{\gamma_{\max}}+\gamma_{\min}}=
\frac{\sqrt{\kappa(\tilde{\mathcal{H}}_h)}-1}{\sqrt{\kappa(\tilde{\mathcal{H}}_h)}+1},
\end{eqnarray*}
where $\kappa(\tilde{\mathcal{H}}_h)$ is the spectral condition number of
$\tilde{\mathcal{H}}_h$.

Now we summarize all the previous results in the form of a theorem.
\begin{theorem}
\label{th-lde-wrhss-time}
Consider the WR-HSS method in time domain
(\ref{lde-wrhss-iteration-1}) for the unsteady discrete elliptic
problem in time domain (\ref{nsa-pd-lde}).
Assume that the data satisfies $f\in\bbb{L}_2^r(\reals)$, and the solution of (\ref{nsa-pd-lde})
satisfies $x_{\star}\in\bbb{L}_2^r(\reals)$, then the solution of (\ref{nsa-pd-lde})
belongs to the Banach space $\bbb{V}_t$, i.e., $x_{\star}\in\bbb{V}_t$.
Let $\alpha$ be a positive constant, and the initial guess of the WR-HSS method
(\ref{lde-wrhss-iteration-1}) belongs to $\bbb{V}_t\cap\bbb{L}_2^r(\reals)$,
then the WR-HSS method (\ref{lde-wrhss-iteration-1})
and its frequency domain counterpart (\ref{lde-wrhss-frequency})
are closed in $\bbb{V}_t\cap\bbb{L}_2^r(\reals)$
and $\bbb{V}_{\omega}\cap\bbb{L}_2^r(\reals)$ respectively.

In addition, the two consecutive iterates of the WR-HSS method in time domain
(\ref{lde-wrhss-iteration-1}) satisfy the following contraction
condition
\begin{eqnarray*}
  \interleave x^{(k+1)}-x_{\star}\interleave_t \le
  \sigma(\alpha)\, \interleave x^{(k)}-x_{\star}\interleave_t,
\end{eqnarray*}
where the upper bound $\sigma(\alpha)$ is given by
\begin{eqnarray*}
\sigma(\alpha) = \max_{\lambda_j\in\lambda(\tilde{\mathcal{H}}_h)}
\left|\frac{\alpha-\lambda_j} {\alpha+\lambda_j}\right|<1,
\end{eqnarray*}
here $\lambda(\tilde{\mathcal{H}}_h)$ is the spectral set of the
matrix $\tilde{\mathcal{H}}_h = H$. Therefore, the WR-HSS method in time domain
(\ref{lde-wrhss-iteration-1}) is convergent for any positive constant
$\alpha$.

Moreover, if $\gamma_{\min}$ and $\gamma_{\max}$ are the lower and
the upper bounds of the eigenvalues of the matrix
$\tilde{\mathcal{H}}_h = H$, respectively, then
\begin{eqnarray*}
\alpha^* = \arg\min_{\alpha}\left\{
\max_{\gamma_{\min}\le\lambda\le\gamma_{\max}}\left|
\frac{\alpha-\lambda_j}{\alpha+\lambda_j} \right|
\right\}=\sqrt{\gamma_{\min}\gamma_{\max}}
\end{eqnarray*}
and
\begin{eqnarray*}
\sigma(\alpha^*) = \frac{\sqrt{\gamma_{\max}}-\gamma_{\min}}
{\sqrt{\gamma_{\max}}+\gamma_{\min}}=
\frac{\sqrt{\kappa(\tilde{\mathcal{H}}_h)}-1}{\sqrt{\kappa(\tilde{\mathcal{H}}_h)}+1},
\end{eqnarray*}
where $\kappa(\tilde{\mathcal{H}}_h)$ is the spectral condition number of
$\tilde{\mathcal{H}}_h$.
\end{theorem}

\section{Implementation details}
\label{sec-implementation-details}

The WR-HSS method discussed in the previous sections is a
continuous-time waveform relaxation method which generate a sequence of
approximate solutions $\{x^{(k)}\}\subset\bbb{L}_2^r(\reals)$
along the whole time axis, i.e., the analytical solution of a
certain system of linear equations and the analytical solution of
a certain system of linear differential equations are required in the
two half steps of each WR-HSS iteration. For
the above reason, the WR-HSS method is therefore mainly of
theoretical interest. In actual implementation, the continuous-time
method should be replaced by a discrete-time method, i.e., the
functions/waveforms are represented discretely as vectors defined on
successive time levels, and the system of linear equations and
the system of linear differential equations are solved by suitable
time-stepping techniques.

We consider the numerical solution of the unsteady elliptic problem
(\ref{uns-elliptic}) on domain $\Omega\subset\reals^d$ and finite
time interval $[0,\mbox{T}]$. The spatial semi-discretization by
using centered difference scheme on equidistant grid with
spatial-step-size $h=\frac{1}{n+1}$ leads to the following unsteady
discrete elliptic problem
\begin{eqnarray}
\label{nsa-pd-lde-finite} \mathcal{L}_h(x) = B\,\dot x(t) + A\, x(t) = f(t),\quad
x(0)=x_0,\quad t\in[0,\mbox{T}]
\end{eqnarray}
with $B=h^2I\in\ceals^{r\times r}$ being Hermitian and $A\in\ceals^{r\times r}$
being non-Hermitian positive definite, here $r=n^d$. The temporal discretization of
the unsteady discrete elliptic problem (\ref{nsa-pd-lde-finite}) by using backward
Euler formula leads to the following difference equations
\begin{eqnarray}
\label{nas-pd-lde-finite-bEuler}
\left(\frac{1}{\Delta t}B+A\right)x_{j+1} - \frac{1}{\Delta t}B\,x_{j} = f_{j+1}, \quad
j = 0,1,2,\ldots,\ell_t,
\end{eqnarray}
where $\Delta t$ is time-step-size, $\ell_t$ is number of time levels,
and $\mbox{T}=\ell_t \times \Delta t$. Moreover, $x_j$ is the approximate value of $x(t)$
on time level $t_j=j\times\Delta t$, and $f_j=f(t_j)$. The above difference equations can be
equivalently rewritten as a discrete linear convolution operator form
\begin{eqnarray}
\label{nas-pd-lde-finite-bEuler-conv}
\mathcal{L}_{\Delta t}\,x_{\Delta t} = f_{\Delta t},
\end{eqnarray}
or equivalently,
\begin{eqnarray*}
\left(\mathcal{L}_{\Delta t}\,x_{\Delta t}\right)_j=\sum_{i=0}^{j}L_{j-i}x_i=f_j,
\quad j=0,1,2,\ldots,\ell_t
\end{eqnarray*}
with matrix-valued kernel
\begin{eqnarray*}
L_{\Delta t} &=& \left\{L_0,L_1,L_2,\ldots,L_{\ell_t}\right\}\\
&=& \left\{\left(\frac{1}{\Delta t}B+A\right),- \frac{1}{\Delta t}B,0,\ldots,0\right\}
\end{eqnarray*}
and vector-valued sequences $x_{\Delta t}=\{x_j\}_{j=0}^{\ell_t}$,
$f_{\Delta t}=\{f_j\}_{j=0}^{\ell_t}$. The discrete linear convolution operator equation
(\ref{nas-pd-lde-finite-bEuler-conv}) can be solved time-level-by-time-level directly.
On each time level, we need to solve only system of linear equations of the form
\begin{eqnarray*}
\left(\frac{h^2}{\Delta t}I+A\right) y_j=c_j
\end{eqnarray*}
by classical band solvers or subspace iterative solvers, e.g.,
GMRES. In general, the cost of classical band solver is $O(n^ds^2)$
where $s$ is the bandwidth and in our context $s=O(n^{d-1})$. When
$d=1$, the above system of linear equations can be solved with
optimal arithmetic by using classical band solver, but it is no
longer true for the cases $d\ge 2$ and $n$ large enough. For these
latter cases, the GMRES is a better choice for the solution of the
above system of linear equations. Specifically, we use the restarted
GMRES($m$) with Householder Arnoldi's procedure, where $m$ is the
restarted parameter. If we have obtained an approximate solution
$y^{(k)}$ of the above system of linear equations by the restarted
GMRES($m$), the corresponding residual vector is define as
\begin{eqnarray*}
r^{(k)} = c_j - \left(\frac{h^2}{\Delta t}I+A\right) y_j^{(k)}.
\end{eqnarray*}
In addition, we denote the solver of the discrete linear convolution
operator equation (\ref{nas-pd-lde-finite-bEuler-conv}) by directly
using the restarted GMRES($m$) on each time level as DGMRES.

As stated in Section \ref{sec-hss}, the WR method is another way to solve the discrete
elliptic problem
(\ref{nsa-pd-lde-finite}), i.e.,
\begin{eqnarray*}
M_B\,\dot x^{(k+1)}(t) + M_A\, x^{(k+1)}(t) = N_B\,\dot x^{(k)}(t) + N_A\, x^{(k)}(t) + f(t),
\end{eqnarray*}
and it can be represented as a one-step operator splitting iterative method.
The temporal discretization of the above one-step WR method by using backward Euler formula
leads to the following discrete-time WR method
\begin{eqnarray*}
\left(\frac{1}{\Delta t}M_B+M_A\right)x_{j+1}^{(k+1)} - \frac{1}{\Delta t}M_B\,x_{j}^{(k+1)}=
\left(\frac{1}{\Delta t}N_B+M_A\right)x_{j+1}^{(k)} - \frac{1}{\Delta t}N_B\,x_{j}^{(k)} + f_{j+1}.
\end{eqnarray*}
This discrete-time WR method can be
equivalently rewritten as a discrete linear convolution operator form
\begin{eqnarray}
\label{dwr-operator-splitting}
\mathcal{M}_{\Delta t}\,x_{\Delta t}^{(k+1)} = \mathcal{N}_{\Delta t}\,x_{\Delta t}^{(k)} + f_{\Delta t}
\end{eqnarray}
with matrix-valued kernel
\begin{eqnarray}
\label{matrix-valued-kernel-M}
M_{\Delta t} &=& \left\{M_0,M_1,M_2,\ldots,M_{\ell_t}\right\} \nonumber \\
&=& \left\{\left(\frac{1}{\Delta t}M_B+M_A\right),- \frac{1}{\Delta t}M_B,0,\ldots,0\right\}
\end{eqnarray}
and
\begin{eqnarray}
\label{matrix-valued-kernel-N}
N_{\Delta t} &=& \left\{N_0,N_1,N_2,\ldots,N_{\ell_t}\right\} \nonumber \\
&=& \left\{\left(\frac{1}{\Delta t}N_B+N_A\right),- \frac{1}{\Delta t}N_B,0,\ldots,0\right\}.
\end{eqnarray}
Obviously, the discrete-time WR method (\ref{dwr-operator-splitting}) is a one-step operator
splitting iterative method for solving the discrete linear convolution operator equation
(\ref{nas-pd-lde-finite-bEuler-conv}). One of the special case of this kind of operator
splitting iterative method is the WR-SOR method, which is based on the matrices splitting
\begin{eqnarray*}
B = h^2I - 0\quad\mbox{and}\quad A = \left(\frac{1}{\tau}D_A-L_A\right)
-\left(\frac{1-\tau}{\tau}D_A+U_A\right),
\end{eqnarray*}
where $\tau$ is an iterative parameter, $D_A$ is diagonal, $L_A$ is
strictly lower triangular, and $U_A$ is strictly upper triangular.
The solution of the discrete linear convolution operator equation
(\ref{dwr-operator-splitting}) is based on the solution of a series of systems of
linear equations
\begin{eqnarray*}
\left(\frac{h^2}{\Delta t}I+\frac{1}{\tau}D_A-L_A\right)\,y_j=c_j
\end{eqnarray*}
along each time level. Obviously, the coefficient matrix of the above system of linear equations is lower triangular,
and only $O(n^ds)$ operations are needed for the solution of it, where $s$ is the bandwidth
of the lower triangular matrix.

We note that the WR-HSS method is a two-step operator splitting iterative method for solving discrete
elliptic problem (\ref{nsa-pd-lde-finite}). Therefore, we describe the general framework of
the two-step operator splitting iterative method and its temporal discretization. Based on the following
two operators splitting of linear differential operator $\mathcal{L}_h$
\begin{eqnarray*}
\mathcal{L}_h &=& \mathcal{M}_1-\mathcal{N}_1 = \left(M_{B_1}\frac{\mbox{d}}{\mbox{d}t}+M_{A_1}\right)-\left(N_{B_1}\frac{\mbox{d}}{\mbox{d}t}+N_{A_1}\right) \\
&=& \mathcal{M}_2-\mathcal{N}_2 = \left(M_{B_2}\frac{\mbox{d}}{\mbox{d}t}+M_{A_2}\right)-\left(N_{B_2}\frac{\mbox{d}}{\mbox{d}t}+N_{A_2}\right)
\end{eqnarray*}
with matrices splitting
\begin{eqnarray*}
\begin{array}{rcl}
B & = & M_{B_1}-N_{B_1} \\
&=& M_{B_2}-N_{B_2}
\end{array}
\quad\mbox{and}\quad
\begin{array}{rcl}
A & = & M_{A_1}-N_{A_1} \\
&=& M_{A_2}-N_{A_2}
\end{array},
\end{eqnarray*}
we can define the continuous-time two-step operator splitting iterative method as
\begin{eqnarray*}
\left\{
\begin{array}{lll}
  \mathcal{M}_1\,x^{(k+\frac{1}{2})} &=& \mathcal{N}_1\,x^{(k)}+f,\\
  \mathcal{M}_2\,x^{(k+1)} &=& \mathcal{N}_2\,x^{(k+\frac{1}{2})}+f.
\end{array}
\right.
\end{eqnarray*}
After temporal discretization of the continuous-time two-step operator splitting iterative method
by using backward Euler formula, we have the discrete-time two-step operator splitting
iterative method
\begin{eqnarray*}
\left\{
\begin{array}{lll}
  \mathcal{M}_{1,\Delta t}\,x_{\Delta t}^{(k+\frac{1}{2})} &=& \mathcal{N}_{1,\Delta t}\,x_{\Delta t}^{(k)}+f_{\Delta t},\\
  \mathcal{M}_{2,\Delta t}\,x_{\Delta t}^{(k+1)} &=& \mathcal{N}_{2,\Delta t}\,x_{\Delta t}^{(k+\frac{1}{2})}+f_{\Delta t},
\end{array}
\right.
\end{eqnarray*}
where $\mathcal{M}_{i,\Delta t}$ and $\mathcal{N}_{i,\Delta t}$, $i=1,2$, are discrete
linear convolution operators with matrix-valued kernels $M_{i,\Delta t}$ and $N_{i,\Delta t}$,
$i=1,2$, defined similarly to (\ref{matrix-valued-kernel-M}) and (\ref{matrix-valued-kernel-N}).
Obviously, the WR-HSS method can be considered as a special case of the above two-step operator
iterative splitting method with matrices splitting
\begin{eqnarray*}
\begin{array}{rcl}
B & = & 0-(-h^2I) \\
&=& h^2I-0
\end{array}
\quad\mbox{and}\quad
\begin{array}{rcl}
A & = & (\alpha\,I+H)-(\alpha\,I-S) \\
&=& (\alpha\,I+S)-(\alpha\,I-H)
\end{array}.
\end{eqnarray*}
During each iteration of the WR-HSS method, the solution of two series of systems of linear equations
along each time level are involved
\begin{eqnarray}
\left(\alpha\,I+H\right)\,y_j=c_j, \label{shifted-hermitian-coe}\\
\left(\left(\frac{h^2}{\Delta t}+\alpha\right)\,I+S\right)\,z_j=b_j. \label{shifted-skew-hermitian-coe}
\end{eqnarray}
Since the matrices $H$ and $S$ are Hermitian part and skew-Hermitian part of the coefficient matrix
$A$ which is arising from the spatial semi-discretization of the unsteady elliptic problem
(\ref{uns-elliptic}) by the centered difference scheme, the systems of
linear equations (\ref{shifted-hermitian-coe}) and
(\ref{shifted-skew-hermitian-coe}) are solved efficiently by the sine and the modified sine
transforms with only $O(n^d\log n)$ operations in our context, respectively. We remark that
we can do better with cyclic reduction (see \cite{Elman90,Elman91}) or multigrid methods
(see \cite{Hackbusch85}) in $O(n^d)$ operations for solving
(\ref{shifted-hermitian-coe}) and (\ref{shifted-skew-hermitian-coe}).

In order to make the WR method more efficient and more practical in actual
implementation, the windowing technique is frequently introduced to the WR
method. Specifically, windowing technique is to divide the whole long time
interval into a number of short time subintervals, and apply the
WR method on each subinterval. For solving the unsteady discrete elliptic problem
(\ref{nsa-pd-lde-finite}) on finite time interval $[0,\mbox{T}]$, we choose $J+1$ time levels,
i.e., $0=\txt{t}_0<\txt{t}_1<\cdots<\txt{t}_J=\txt{T}$, to divide time interval $[0,\mbox{T}]$
into $J$ smaller equidistance subintervals $(\txt{t}_{i-1},\txt{t}_{i}]$, $i=1,2,\ldots,J$,
with $\ell_{t,J}$ time levels on each subinterval $(\txt{t}_{i-1},\txt{t}_{i}]$ and
$\ell_{t,J} \times J = \ell_t$, then the WR method, such as the WR-SOR method and
the WR-HSS method, can be applied to solve the unsteady discrete elliptic problem (\ref{nsa-pd-lde-finite})
on each subinterval $(\txt{t}_{i-1},\txt{t}_{i}]$.
Since the subintervals are shorter, fewer number of time levels are involved,
the number of iterations of the WR method applied on each subinterval is smaller than
that of the WR method applied on the whole long time interval. Therefore, the overall computation loads on all of the
subintervals is smaller than the computation loads while simulating once and for all
on the whole long time interval.

\section{Numerical examples}
\label{sec-num-example}

In this section, we present some numerical examples to demonstrate
the correctness of the previously proposed theory and the
effectiveness of the WR-HSS method.

If $z_{\Delta
t}=\{z_j\}_{j=0}^{\ell}$ is a vector sequence of length $\ell+1$
with $z_j\in\ceals^r$, the norm of this vector sequence $z_{\Delta
t}$ is defined as
\begin{eqnarray*}
\|z_{\Delta t}\|_p = \left\{
\begin{array}{ll}
    \root \of {\sum _{i=0}^{\ell}\|z_i\|_2^2} & p = 2,\\
    \sup \limits_{0\le i <\ell}\{\|z_i\|_{\infty}\} & p = \infty.
\end{array}
\right.
\end{eqnarray*}

Suppose that we have obtained an approximate solution of the
discrete linear convolution operator equation
(\ref{nas-pd-lde-finite-bEuler-conv}) by some discrete-time WR
method, say $x_{\Delta t}^{(k)}=\{x_j^{(k)}\}_{j=0}^{\ell_t}$, we
define the relative error of the approximate solution $x_{\Delta
t}^{(k)}$ as
\begin{eqnarray*}
\mbox{ERR}=\frac{\|x_{\Delta t}^{(k)}-x_{\Delta t}\|_{\infty}}{\|x_{\Delta t}\|_{\infty}}.
\end{eqnarray*}
In addition, we define the residual vector-valued sequence of the
discrete linear convolution operator equation
(\ref{nas-pd-lde-finite-bEuler-conv}) with respect to $x_{\Delta
t}^{(k)}$ as
\begin{eqnarray*}
r_{\Delta t}^{(k)} = f_{\Delta t} - \mathcal{L}_{\Delta t}\,x_{\Delta t}^{(k)},
\end{eqnarray*}
and the corresponding relative residual is defined as
\begin{eqnarray*}
\mbox{RES}=\frac{\|r_{\Delta t}^{(k)}\|_2}{\|r_{\Delta t}^{(0)}\|_2}.
\end{eqnarray*}

All computations were completed with MATLAB 2014a installed in Windows XP Professional 2002
Service Pack 3 on Intel(R) Core(TM) i3-2130 CPU @ 3.40GHz 3.39GHz with 3.35GB RAM.

\subsection{The 1-dimensional case}
\label{sec-num-1D}

In this subsection, we consider the 1-dimensional unsteady elliptic problem
\begin{eqnarray*}
\frac{\partial u(x,t)}{\partial t} - \frac{\partial ^2 u(x,t)}{\partial x^2}
+ q\, \frac{\partial u(x,t)}{\partial x} = 0
\end{eqnarray*}
on spatial domain $\Omega = [0,1]$ and time interval $[0,\mbox{T}]$,
with constant coefficient $q$ of the convection term, and subject to
Dirichlet boundary condition. When the centered difference scheme is
applied to the above unsteady elliptic problem, and the natural lexicographic
ordering is employed to the unknowns, we get the unsteady
discrete elliptic problem with coefficients
\begin{eqnarray*}
B=h^2\, I\in\reals^{n\times n} \quad\mbox{and}\quad
A=\tdiag (-1-\mbox{\it Re}, 2, -1+\mbox{\it Re})\in\reals^{n\times n},
\end{eqnarray*}
where $\mbox{\it Re}=\frac{qh}{2}$ is the mesh Reynolds number. For the convenience of
error comparison, the exact solution of the corresponding unsteady discrete
elliptic problem is artificially chosen to be
\begin{eqnarray*}
x(t) = e^{-t}\,\mbox{\bf 1},
\end{eqnarray*}
where $\mbox{\bf 1} = (1,1,\ldots, 1)^T$.
In the tests, we choose
$n=64$ and $\mbox{T}=2$. According to Remark \ref{rema-wrhss-convergence}, we have
\begin{eqnarray*}
\rho(\mathcal{K}_{\mbox{\tiny WR-HSS}}) = \sup_{\omega\in\reals}
    \rho(\tilde K_{\mbox{\tiny WR-HSS}})
    \approx \sup_{-\omega_c\le\omega\le \omega_c}
    \rho(\tilde K_{\mbox{\tiny WR-HSS}})
    \le \sigma(\alpha),
\end{eqnarray*}
where $\omega_c$ is a given large upper bound of frequency $\omega$.
The above fact demonstrates that the values of spectral radius
$\rho(\tilde K_{\mbox{\tiny WR-HSS}})$ in frequency domain
represent the value of spectral radius $\rho(\mathcal{K}_{\mbox{\tiny WR-HSS}})$
in time domain to some extent.

\begin{figure}
            \centering
            \begin{tabular}{c}
                \includegraphics[scale=0.96,bb=100 350 520 520]{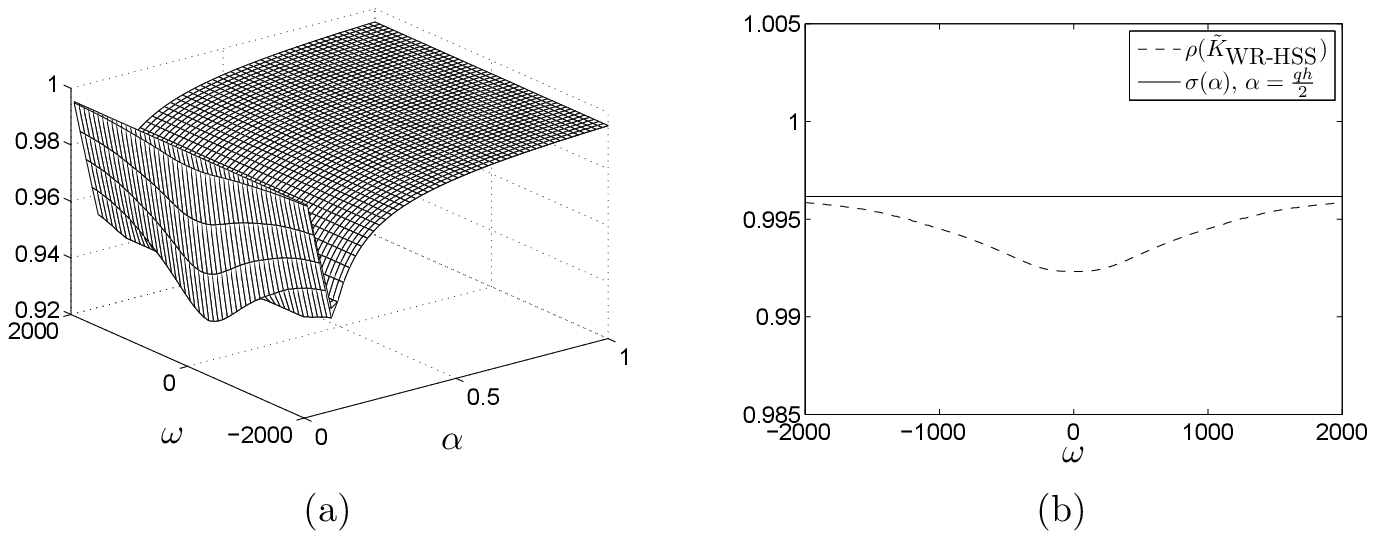}
                \\
            \end{tabular}
            \caption{For $q=1$: (a) surfaces of the spectral radius $\rho(\tilde K_{\mbox{\tiny WR-HSS}})$
            and the upper bound $\sigma(\alpha)$; (b) sectional drawing of the previous
            surfaces for $\alpha = \frac{qh}{2}$, dashed line
            for $\rho(\tilde K_{\mbox{\tiny WR-HSS}})$ and solid line for $\sigma(\alpha)$.
            }
            \label{fig:radii-frequency-centered64-q1}
\end{figure}

\begin{figure}
            \centering
            \begin{tabular}{c}
                \includegraphics[scale=0.96,bb=100 350 520 520]{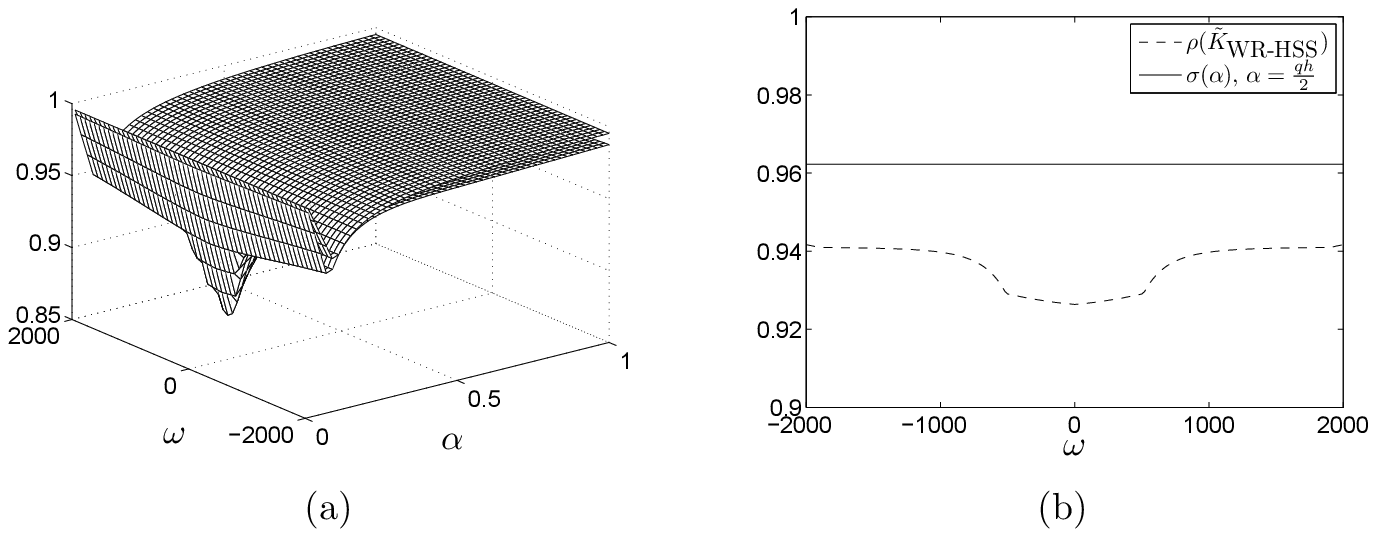}
                \\
            \end{tabular}
            \caption{For $q=10$: (a) surfaces of the spectral radius $\rho(\tilde K_{\mbox{\tiny WR-HSS}})$
            and the upper bound $\sigma(\alpha)$; (b) sectional drawing of the previous
            surfaces for $\alpha = \frac{qh}{2}$, dashed line
            for $\rho(\tilde K_{\mbox{\tiny WR-HSS}})$ and solid line for $\sigma(\alpha)$.
            }
            \label{fig:radii-frequency-centered64-q10}
\end{figure}

\begin{figure}
            \centering
            \begin{tabular}{c}
                \includegraphics[scale=0.96,bb=100 350 520 520]{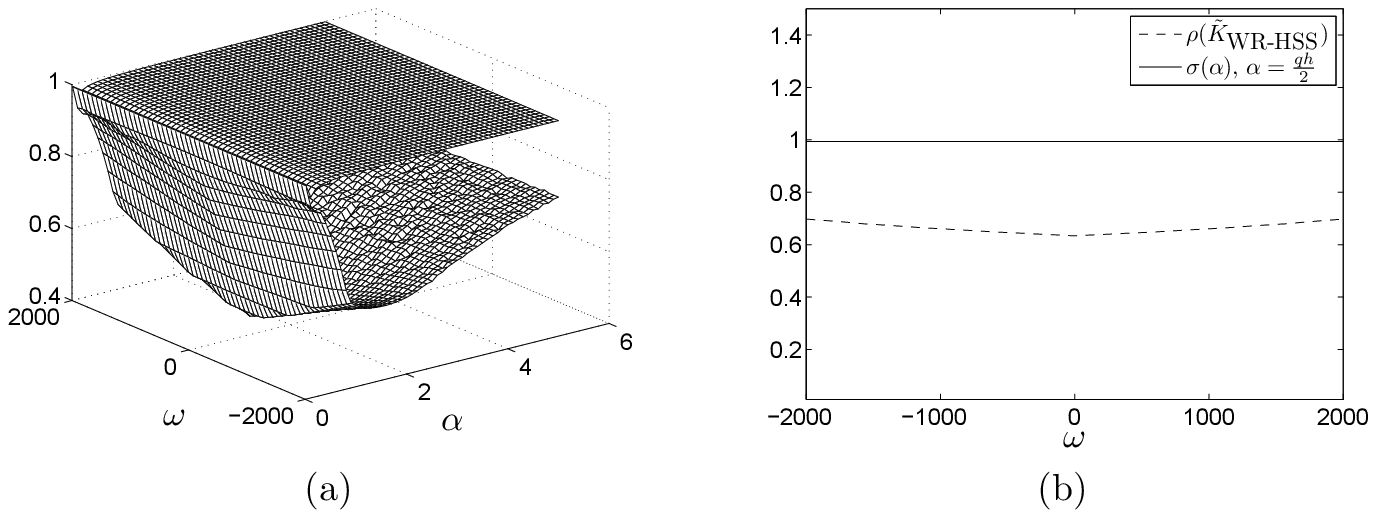}
                \\
            \end{tabular}
            \caption{For $q=100$: (a) surfaces of the spectral radius $\rho(\tilde K_{\mbox{\tiny WR-HSS}})$
            and the upper bound $\sigma(\alpha)$; (b) sectional drawing of the previous
            surfaces for $\alpha = \frac{qh}{2}$, dashed line
            for $\rho(\tilde K_{\mbox{\tiny WR-HSS}})$ and solid line for $\sigma(\alpha)$.
            }
            \label{fig:radii-frequency-centered64-q100}
\end{figure}

\begin{figure}
            \centering
            \begin{tabular}{c}
                \includegraphics[scale=0.96,bb=100 350 520 520]{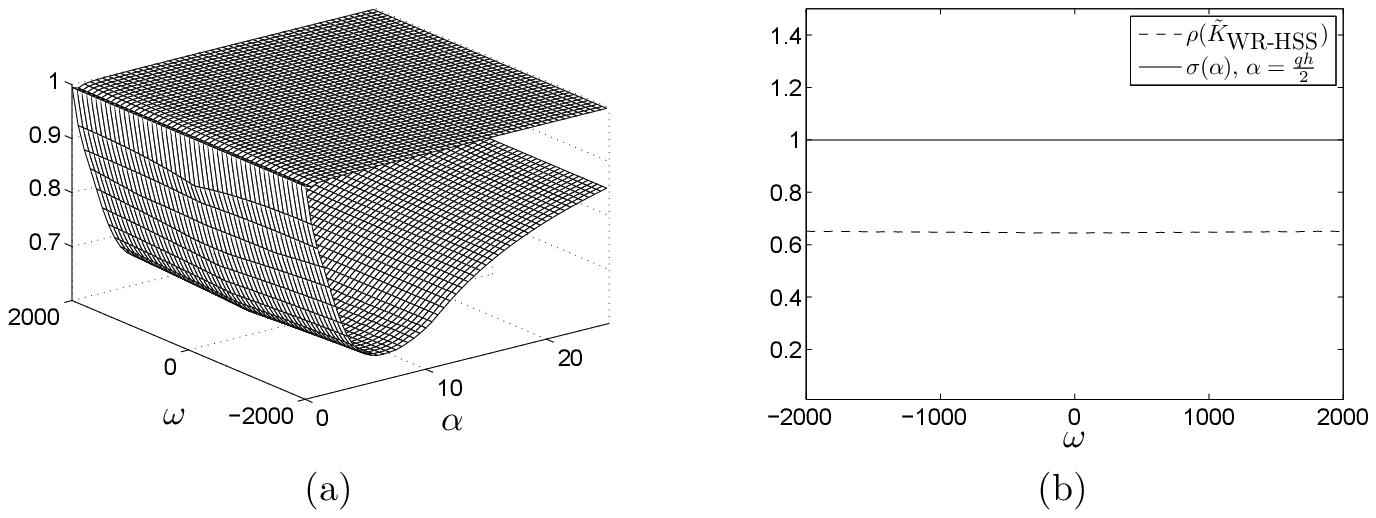}
                \\
            \end{tabular}
            \caption{For $q=1000$: (a) surfaces of the spectral radius $\rho(\tilde K_{\mbox{\tiny WR-HSS}})$
            and the upper bound $\sigma(\alpha)$; (b) sectional drawing of the previous
            surfaces for $\alpha = \frac{qh}{2}$, dashed line
            for $\rho(\tilde K_{\mbox{\tiny WR-HSS}})$ and solid line for $\sigma(\alpha)$.
            }
            \label{fig:radii-frequency-centered64-q1000}
\end{figure}

Figures
\ref{fig:radii-frequency-centered64-q1}-\ref{fig:radii-frequency-centered64-q1000}
show the surfaces of the spectral
radius $\rho(\tilde K_{\mbox{\tiny WR-HSS}})$ and the upper bound $\sigma(\alpha)$
on $\omega$-$\alpha$-plane with $\omega_c = 2000$ for different values of $q$,
and the corresponding sectional drawing of the previous surfaces for
$\alpha = \frac{qh}{2}$. In addition, the interval of $\alpha$ is determined
accordingly. When $q$ is small (e.g., $q=1$ in
Figure \ref{fig:radii-frequency-centered64-q1}), the surfaces in sub-figure-(a)
stick together, and it is difficult to tell the difference between them. Moreover,
the corresponding sectional drawing in sub-figure-(b) gives a better illustration
of the tiny difference between the spectral
radius $\rho(\tilde K_{\mbox{\tiny WR-HSS}})$ and the upper bound $\sigma(\alpha)$.
When $q$ becomes larger (e.g., $q=1000$ in
Figure \ref{fig:radii-frequency-centered64-q1000}), the surfaces in sub-figure-(a)
are wide apart from each other, and the corresponding sectional drawing in
sub-figure-(b) also demonstrates a larger difference between the spectral
radius $\rho(\tilde K_{\mbox{\tiny WR-HSS}})$ and the upper bound $\sigma(\alpha)$.

Table \ref{tab:values_of_sigma_and_radii} lists the values of the upper bound
$\sigma(\alpha)$ and the intervals of the spectral radius
$\rho(\tilde K_{\mbox{\tiny WR-HSS}})$ with $-\omega_c\le\omega\le \omega_c$
for $\alpha=\frac{qh}{2}$ and different
values of $q$. Obviously, the values of the upper bound $\sigma(\alpha)$ are all
less than but close to one, which means that the convergence of the WR-HSS method
is guaranteed, but the actual convergence rate can not be revealed correctly. In addition,
the value of the spectral radius $\rho(\tilde K_{\mbox{\tiny WR-HSS}})$ is close
to the value of the upper bound $\sigma(\alpha)$ for small $q$, and
the value of the spectral radius $\rho(\tilde K_{\mbox{\tiny WR-HSS}})$ decreases
fast, when the value of $q$ increases.

\begin{figure}
            \centering
            \begin{tabular}{c}
                \includegraphics[scale=0.6]{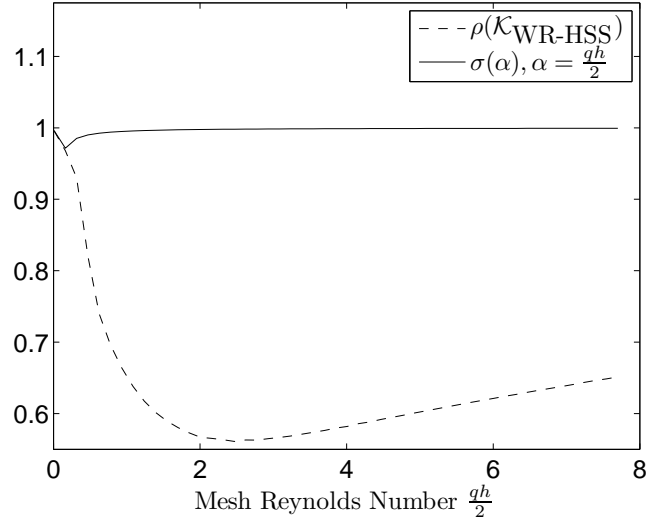}
                \\
            \end{tabular}
            \caption{The spectral radius
            $\rho(\mathcal{K}_{\mbox{\tiny WR-HSS}})$
            and the upper bound $\sigma(\alpha)$ with respect to
            the mesh Reynolds number $\frac{qh}{2}$, dashed line
            for $\rho(\mathcal{K}_{\mbox{\tiny WR-HSS}})$ and solid line for $\sigma(\alpha)$.
            }
            \label{fig:radii-upperboundq-centered64}
\end{figure}

Figure \ref{fig:radii-upperboundq-centered64} depict the curves of the spectral radius
$\rho(\mathcal{K}_{\mbox{\tiny WR-HSS}})$ and the upper bound $\sigma(\alpha)$ with
respect to the mesh Reynolds number $\frac{qh}{2}$ with $\alpha = \frac{qh}{2}$.
We find that the two curves stay close when the mesh Reynolds number $\frac{qh}{2}$
is small, and when the mesh Reynolds number $\frac{qh}{2}$ increases, the two curves
are apart from each other rapidly. Moreover, the curve of the spectral radius
$\rho(\mathcal{K}_{\mbox{\tiny WR-HSS}})$ stays below the curve of
the upper bound $\sigma(\alpha)$ all the time.

The above observations show that the convergence of the WR-HSS
method is unconditionally guaranteed for any positive parameter
$\alpha$, the upper bound $\sigma(\alpha)$ is close to the spectral
radius $\rho(\mathcal{K}_{\mbox{\tiny WR-HSS}})$ for small $q$, and
they are all close to one, which means that $\sigma(\alpha)$ is a
good approximation of the spectral radius $\rho(\mathcal{K}_{\mbox{\tiny WR-HSS}})$ when
the unsteady elliptic problem has a weak convection term, however,
the convergence rate of the WR-HSS method is very slow in this case. When $q$
becomes larger, or say the unsteady elliptic problem has a stronger
convection term, the upper bound $\sigma(\alpha)$ keeps close to
one, but the spectral radius $\rho(\mathcal{K}_{\mbox{\tiny
WR-HSS}})$ is far less than one, which means that the convergence rate
of the WR-HSS method is much faster than the upper bound
$\sigma(\alpha)$ can reveal.

\subsection{The 2-dimensional case}
\label{sec-num-2D}

In this subsection, we compare the WR-HSS method with the DGMRES and
the WR-SOR method to demonstrate the robustness of the WR-HSS
method. We consider the 2-dimensional unsteady elliptic problem
\begin{eqnarray*}
\frac{\partial u(x,t)}{\partial t}-\nabla \cdot [a(x,t)\nabla u(x,t)] +
  \sum _{j=1}^{d} \frac{\partial}{\partial x_j}(q(x,t)u(x,t)) = 0
\end{eqnarray*}
on spatial domain $\Omega$ and time interval $[0,\mbox{T}]$,
with positive constant function $a(x,t)=1$ and constant Reynolds function $q(x,t)=q$,
and subject to Dirichlet boundary condition. The spatial domain $\Omega$ can be
chosen to be a square domain $Q=[0,1]\times[0,1]$ or L-shaped domain $L=Q\setminus G$,
with $G=[0,0.5]\times[0,0.5]$. We remark here that we did not observe any difference
in the quality of numerical results on the square domain and the L-shaped domain.
Thus, we only report the numerical results on square domain.

According to the analysis in Subsection \ref{sec-num-1D},
the WR-HSS method converges fast for solving the unsteady elliptic problem with strong
convection term, thus, $q$ is chosen to be large in the tests in this subsection.
When the centered difference scheme is
applied to the above unsteady elliptic problem, and the natural lexicographic
ordering is employed to the unknowns, we get the unsteady
discrete elliptic problem with coefficients
\begin{eqnarray*}
B=h^2\, I\otimes I \quad\mbox{and}\quad
A=I\otimes T_n+ T_n \otimes I,
\end{eqnarray*}
where $T_n=\tdiag (-1-\mbox{\it Re}, 2, -1+\mbox{\it Re})$ with $\mbox{\it Re}=\frac{qh}{2}$
as the mesh Reynolds number. For the convenience of
error comparison, the exact solution of the corresponding unsteady discrete
elliptic problem is artificially chosen to be
\begin{eqnarray*}
x = e^{-t}\,\mbox{\bf 1},
\end{eqnarray*}
where $\mbox{\bf 1} = (1,1,\ldots, 1)^T$.

In the tests, we choose spatial-grid-size $n=127,255,511,1023$ and time
interval $[0,\mbox{T}]=[0,\ell_t\times\Delta t]$. The system size on
each time level varies from $O(10^{4})$ to $O(10^6)$. The stopping
criterion of the discrete-time WR method on each window is given by
\begin{eqnarray*}
\frac{\|r_{\Delta t}^{(k)}\|_2}{\|r_{\Delta
t}^{(0)}\|_2}<\varepsilon,
\end{eqnarray*}
where $\varepsilon$ is a tolerance to control the above stopping
criterion. In addition, The stopping criterion of the restarted
GMRES($m$) used on each time level of the DGMRES is given by
\begin{eqnarray*}
\frac{\|r^{(k)}\|_2}{\|r^{(0)}\|_2}<\eta,
\end{eqnarray*}
where $\eta$ is a tolerance to control the above stopping
criterion.

Tables
\ref{tab:interval_of_parameter_dtMinus4}-\ref{tab:interval_of_parameter_dtMinus6}
list the experimental feasible interval of iteration parameters for
the WR-SOR method and the WR-HSS method for the settings $\Delta t
=10^{-4},10^{-5},10^{-6}$, $\ell_t=\ell_{t,J} \times J=5\times 1$,
$q=2000,3000$ and $n=127,255$. Since these are feasibility tests, the
corresponding tolerance is set to be $\varepsilon=0.05$. Here in the tables,
``$\tau$'' is the parameter of the WR-SOR method, and ``$\alpha$''
is that of the WR-HSS method. Moreover, ``$100+$'' means that
the tested parameter $\alpha$ exceeds $100$. Obviously, the length
of the feasible interval of $\alpha$ is much larger than the length
of the feasible interval of $\tau$ in all cases. This means that
the iterative parameter of the
WR-SOR method is much more sensitive than that of the WR-HSS
method. Thus, the WR-HSS method is a better choice for practical
aspect. In addition, we observe that the WR-HSS method is not
feasible for all positive $\alpha$ in the tests, which contradicts
with the result stated in Theorem \ref{th-lde-wrhss-frequency} and
the observation in Subsection \ref{sec-num-1D}. The
possible reason is that Theorem \ref{th-lde-wrhss-frequency}
describes the convergence behavior of the continuous-time WR-HSS
method, but the actual implemented method is the discrete-time
WR-HSS method which is not a HS splitting based method in essential,
therefore, not all positive $\alpha$ leads to a convergent
discrete-time WR-HSS method.

Tables
\ref{tab:all_wr_methods_table_dtMinus4_q2000_Lt25}-\ref{tab:all_wr_methods_table_dtMinus6_q3000_Lt25}
list the numerical results of the WR-HSS method, the WR-SOR method
and the DGMRES. The settings of the problems are given by $\Delta t
=10^{-4},10^{-5},10^{-6}$, $\ell_t=\ell_{t,J} \times J=5\times 5$,
$q=2000,3000$ and $n=255,511,1023$. The iterative parameter
$\alpha$ of the WR-HSS method is given by $\alpha=\frac{qh}{2}$, the
iterative parameter $\tau$ of the WR-SOR method is determined
experimentally, and the restarted parameter $m$ used in the DGMRES
on each time level is set to be $m=5$ which is the largest $m$ we can
use to avoid running out of memory in all cases. In these tables, ``IT''
represents the average number of iterations on each window of the
WR-HSS method and the WR-SOR method, and the average number of
iterations (i.e., the average number of matrix-vector product)
on each time level of the DGMRES. The maximum number of
iterations is $7000$, and ``CPU'' represents the total computation
time on the whole time interval. In addition, the tolerance of the
WR-HSS method and the WR-SOR method on each window is fixed to be
$\varepsilon=10^{-5}$, and the tolerance of the DGMRES on each time
level is given by
\begin{eqnarray*}
\eta=\left\{
\begin{array}{lll}
10^{-8}, & \mbox{for} & \Delta t = 10^{-4}, \\
10^{-9}, & \mbox{for} & \Delta t = 10^{-5}, \\
10^{-10}, & \mbox{for} & \Delta t = 10^{-6}, \\
\end{array}
\right.
\end{eqnarray*}
such that the relative error ``ERR'' of the approximate solution
obtained by the DGMRES can be comparable to that of the approximate
solution obtained by the WR-HSS method.

According to Tables
\ref{tab:all_wr_methods_table_dtMinus4_q2000_Lt25}-\ref{tab:all_wr_methods_table_dtMinus6_q3000_Lt25},
we find a general fact that the WR-HSS method and the DGMRES
outperform the WR-SOR method in the aspects of both the IT and the
CPU in all cases. In most cases, the WR-SOR method did not attain
the given tolerance $\varepsilon=10^{-5}$ after the maximum number
of iterations, i.e., $\mbox{IT}=7000$. In the sequel, we only discuss the
numerical behavior of the WR-HSS method and the DGMRES.

For $\Delta t=10^{-4}, 10^{-5}$ and $q=2000,3000$, the IT of the
DGMRES is at least 6 times as many as the IT of the WR-HSS method,
and the CPU of the DGMRES is twice the CPU of the WR-HSS method.
Obviously, the WR-HSS method is much more efficient than the DGMRES
in these cases.

For $\Delta t=10^{-6}$, $q=2000$, $n=255$, both the IT and the CPU
of the WR-HSS method are larger than the IT and the CPU of the
DGMRES. For $\Delta t=10^{-6}$, $q=2000$, $n=511$ and $\Delta
t=10^{-6}$, $q=3000$, $n=255$, the IT of the WR-HSS method is less
than the IT of the DGMRES, but the CPU of the WR-HSS method is
larger than the CPU of the DGMRES. It seems that the DGMRES is more
efficient than the WR-HSS method in these cases, but the fact behind
the illusion is that the ERR of the DGMRES is about 100 times as
large as the ERR of the WR-HSS method. If we decrease the tolerance
$\eta$ of the DGMRES on each time level such that the same ERR as
the WR-HSS method attained for the DGMRES in these cases, it can
make all the difference.

In all cases in Tables
\ref{tab:all_wr_methods_table_dtMinus4_q2000_Lt25}-\ref{tab:all_wr_methods_table_dtMinus6_q3000_Lt25},
the IT and the CPU of the DGMRES both increase proportionally to the
increasing proportion of the spatial-grid-size, and they slightly
increase with $q$. For $\Delta t=10^{-4}, 10^{-5}$ and
$q=2000,3000$, the increasing behavior of the WR-HSS method is
similar to that of the DGMRES. However, for $\Delta t=10^{-6}$ and
$q=2000,3000$, the surprising fact is that the IT of the WR-HSS method
is slightly decreasing with the increasing spatial-grid-size.

Tables
\ref{tab:long_time_interval_dtMinus3}-\ref{tab:long_time_interval_dtMinus5}
list the results of the numerical solution of the 2-dimensional unsteady
elliptic problem on the long time interval by using the WR-HSS method
for the settings $[0,\mbox{T}]=[0,1]$, $\Delta
t=10^{-3},10^{-4},10^{-5}$, $\ell_t=\ell_{t,J} \times J=5\times 200,
5\times 2000, 5\times 20000$, $q=2000,3000$, and $n=127,255$. The
iterative parameter $\alpha$ of the WR-HSS method is given by
$\alpha=\frac{qh}{2}$, and the tolerance of the WR-HSS method on
each window is given by $\varepsilon=10^{-5}$. In these tables, we
find that the IT of the WR-HSS method increases both with $q$ and
with the spatial-grid-size $n+1$. However, the IT of the WR-HSS method
decreases fast with the time-step-size $\Delta t$. Moreover, the ERR
of the WR-HSS method on long time interval remains the same orders
of magnitude as the ERR of the WR-HSS method on short time interval,
e.g., in the cases $\Delta t=10^{-4},10^{-5}$, $q=2000,3000$
and $n=255$.

\section{Conclusions}
\label{conclusions}

In this paper, we propose a two-step operator splitting
method based on the HS splitting of the differential operator for
solving the unsteady discrete elliptic problem, i.e., WR-HSS method.
The theoretical analysis and the numerical results both suggest that the WR-HSS
method is effective for handling the unsteady discrete elliptic
problem. 

The advantages of the WR-HSS method can be addressed in two main aspects.
Firstly, compared with other analytical methods (such as the WR-SOR method),
the WR-HSS method converges unconditionally to the solution of the system 
of linear differential equations, and there is an easy implemented strategy
for computing the iterative parameter. For the WR-SOR method, however,
the length of the feasible interval of the iterative parameter is small,
and the practical iterative parameter is hard to determine.
Secondly, compared with the classical time-stepping methods (such as the DGMRES method) 
who deal with the original system of linear differential equations directly, 
the WR-HSS method splits the original difficult system of linear differential
equations into two easy sub-systems which are much easier to resolve.

In practical aspect, the WR-HSS method must be implemented discretely,
i.e., the discrete-time WR-HSS method, but the theoretical analysis
in this paper are only concentrated on the continuous-time WR-HSS method,
therefore, our future work should be the discussion on the convergence
property of the discrete-time method, and the relationship
between the discrete-time and continuous-time method.



\clearpage


\clearpage

\begin{table}
\setlength{\abovecaptionskip}{0pt}
\setlength{\belowcaptionskip}{10pt}
\centering{
\caption{\label{tab:values_of_sigma_and_radii} Values of the upper bound $\sigma(\alpha)$ and
the spectral radius $\rho(\tilde K_{\mbox{\tiny WR-HSS}})$ with $-\omega_c\le\omega\le\omega_c$ for $\alpha=\frac{qh}{2}$.}
\begin{tabular}{|l|c|c|}\hline
 & $\sigma(\alpha)$ & $\rho(\tilde K_{\mbox{\tiny WR-HSS}})$\\ \hline
 $q=1$ & 0.9962 & (0.9923, 0.9958) \\ \hline
 $q=10$ & 0.9622 & (0.9264, 0.9417) \\ \hline
 $q=100$ & 0.9939 & (0.6339, 0.6976) \\ \hline
 $q=1000$ & 0.9994 & (0.6444, 0.6515) \\ \hline
\end{tabular}}
\end{table}



\begin{table}
\setlength{\abovecaptionskip}{0pt}
\setlength{\belowcaptionskip}{10pt}
\centering{
\caption{\label{tab:interval_of_parameter_dtMinus4} The interval of feasible iteration parameter of the WR-SOR method and the WR-HSS method: $\Delta t=10^{-4}$. }
\begin{tabular}{|l|ll|ll|}\hline
         & \multicolumn{2}{|c|}{$q=2000$} & \multicolumn{2}{|c|}{$q=3000$}  \\ \hline
         & $r=n^2,n=127$& $r=n^2,n=255$ & $r=n^2,n=127$ & $r=n^2,n=255$ \\ \hline
WR-SOR, $\tau$ &   (0.0019, 0.2294) &  (0.0019, 0.4097) &  (0.0019, 0.1582) &  (0.0019, 0.2924) \\ \hline
WR-HSS, $\alpha$ &   (0.9766, 100+) &    (0.9766, 100+) &    (1.4648, 100+) &    (0.7324, 100+) \\ \hline
\end{tabular}}
\end{table}


\begin{table}
\setlength{\abovecaptionskip}{0pt}
\setlength{\belowcaptionskip}{10pt}
\centering{
\caption{\label{tab:interval_of_parameter_dtMinus5} The interval of feasible iteration parameter of the WR-SOR method and the WR-HSS method: $\Delta t=10^{-5}$. }
\begin{tabular}{|l|ll|ll|}\hline
         & \multicolumn{2}{|c|}{$q=2000$} & \multicolumn{2}{|c|}{$q=3000$}  \\ \hline
         & $r=n^2,n=127$& $r=n^2,n=255$ & $r=n^2,n=127$ & $r=n^2,n=255$ \\ \hline
WR-SOR, $\tau$ &   (0.0019, 0.2540) &  (0.0019, 0.4287) &  (0.0019, 0.1707) &  (0.0019, 0.3024) \\ \hline
WR-HSS, $\alpha$ &   (0.9766, 100+) &    (0.9766, 100+) &    (0.7324, 100+) &    (0.7324, 100+) \\ \hline
\end{tabular}}
\end{table}


\begin{table}
\setlength{\abovecaptionskip}{0pt}
\setlength{\belowcaptionskip}{10pt}
\centering{
\caption{\label{tab:interval_of_parameter_dtMinus6} The interval of feasible iteration parameter of the WR-SOR method and the WR-HSS method: $\Delta t=10^{-6}$. }
\begin{tabular}{|l|ll|ll|}\hline
         & \multicolumn{2}{|c|}{$q=2000$} & \multicolumn{2}{|c|}{$q=3000$}  \\ \hline
         & $r=n^2,n=127$& $r=n^2,n=255$ & $r=n^2,n=127$ & $r=n^2,n=255$ \\ \hline
WR-SOR, $\tau$ &   (0.0019, 2.6814) &  (0.0019, 0.7862) &  (0.0019, 0.6169) &  (0.0019, 0.4448) \\ \hline
WR-HSS, $\alpha$ &   (5.1270, 100+) &    (1.9531, 100+) &    (5.1270, 100+) &    (1.4648, 100+) \\ \hline
\end{tabular}}
\end{table}


\begin{table}
\footnotesize \setlength{\abovecaptionskip}{0pt}
\setlength{\belowcaptionskip}{10pt} \centering{
\caption{\label{tab:all_wr_methods_table_dtMinus4_q2000_Lt25} The
number of iterations, computation time and the relative errors of
the approximate solutions for different methods on time interval
$[0,\ell_{t,J}\times J\times\Delta t]=[0,5\times 5\times\Delta t]$:
$\Delta t=10^{-4}$ and $q=2000$.}
\begin{tabular}{|l|ccc|ccc|ccc|}\hline
 & \multicolumn{3}{|c|}{$r=n^2,n=255$} & \multicolumn{3}{|c|}{$r=n^2,n=511$} & \multicolumn{3}{|c|}{$r=n^2,n=1023$} \\ \hline
 & IT & CPU & ERR & IT & CPU & ERR & IT & CPU & ERR \\ \hline
 WR-HSS & 106 & 188.7031 & 6.3E-08 & 155 & 1348.3750 & 7.0E-08 & 272 & 9727.0313 & 1.2E-07 \\ \hline
 WR-SOR & 7000 & 1254.5313 & 5.1E-04 & 7000 & 4878.5469 & 3.7E-03 & 7000 & 18436.4375 & 2.1E-03 \\ \hline
 DGMRES & 649 & 337.9844 & 1.3E-07 & 1096 & 2490.8594 & 2.3E-07 & 2107 & 19215.8281 & 1.7E-07 \\ \hline
\end{tabular}}
\end{table}


\begin{table}
\footnotesize \setlength{\abovecaptionskip}{0pt}
\setlength{\belowcaptionskip}{10pt} \centering{
\caption{\label{tab:all_wr_methods_table_dtMinus4_q3000_Lt25} The
number of iterations, computation time and the relative errors of
the approximate solutions for different methods on time interval
$[0,\ell_{t,J}\times J\times\Delta t]=[0,5\times 5\times\Delta t]$:
$\Delta t=10^{-4}$ and $q=3000$.}
\begin{tabular}{|l|ccc|ccc|ccc|}\hline
  & \multicolumn{3}{|c|}{$r=n^2,n=255$} & \multicolumn{3}{|c|}{$r=n^2,n=511$} & \multicolumn{3}{|c|}{$r=n^2,n=1023$} \\ \hline
 & IT & CPU & ERR & IT & CPU & ERR & IT & CPU & ERR \\ \hline
 WR-HSS & 126 & 226.6406 & 5.8E-08 & 171 & 1483.0000 & 8.5E-08 & 281 & 10086.8125 & 6.7E-08 \\ \hline
 WR-SOR & 7000 & 1310.3281 & 3.1E-04 & 7000 & 4878.7031 & 3.5E-03 & 7000 & 18430.7344 & 3.4E-03 \\ \hline
 DGMRES & 712 & 372.2188 & 2.7E-07 & 1192 & 2711.0938 & 1.1E-07 & 2149 & 19590.9063 & 3.3E-07 \\ \hline
\end{tabular}}
\end{table}


\begin{table}
\footnotesize \setlength{\abovecaptionskip}{0pt}
\setlength{\belowcaptionskip}{10pt} \centering{
\caption{\label{tab:all_wr_methods_table_dtMinus5_q2000_Lt25} The
number of iterations, computation time and the relative errors of
the approximate solutions for different methods on time interval
$[0,\ell_{t,J}\times J\times\Delta t]=[0,5\times 5\times\Delta t]$:
$\Delta t=10^{-5}$ and $q=2000$.}
\begin{tabular}{|l|ccc|ccc|ccc|ccc|}\hline
  & \multicolumn{3}{|c|}{$r=n^2,n=255$} & \multicolumn{3}{|c|}{$r=n^2,n=511$} & \multicolumn{3}{|c|}{$r=n^2,n=1023$} \\ \hline
 & IT & CPU & ERR & IT & CPU & ERR & IT & CPU & ERR \\ \hline
 WR-HSS & 30 & 54.8281 & 3.2E-09 & 47 & 406.3438 & 8.1E-09 & 91 & 3263.3438 & 7.4E-09 \\ \hline
 WR-SOR & 4037 & 755.7500 & 1.1E-08 & 7000 & 4877.6406 & 1.4E-05 & 7000 & 18427.1875 & 4.7E-06 \\ \hline
 DGMRES & 208 & 107.0313 & 4.7E-08 & 410 & 927.1875 & 6.9E-08 & 822 & 7466.7813 & 5.2E-08 \\ \hline
\end{tabular}}
\end{table}


\begin{table}
\footnotesize \setlength{\abovecaptionskip}{0pt}
\setlength{\belowcaptionskip}{10pt} \centering{
\caption{\label{tab:all_wr_methods_table_dtMinus5_q3000_Lt25} The
number of iterations, computation time and the relative errors of
the approximate solutions for different methods on time interval
$[0,\ell_{t,J}\times J\times\Delta t]=[0,5\times 5\times\Delta t]$:
$\Delta t=10^{-5}$ and $q=3000$.}
\begin{tabular}{|l|ccc|ccc|ccc|}\hline
  & \multicolumn{3}{|c|}{$r=n^2,n=255$} & \multicolumn{3}{|c|}{$r=n^2,n=511$} & \multicolumn{3}{|c|}{$r=n^2,n=1023$} \\ \hline
 & IT & CPU & ERR & IT & CPU & ERR & IT & CPU & ERR \\ \hline
 WR-HSS & 45 & 80.8125 & 3.7E-09 & 69 & 592.3594 & 1.2E-08 & 132 & 4736.0313 & 1.2E-08 \\ \hline
 WR-SOR & 3858 & 722.9531 & 1.4E-08 & 7000 & 4886.3906 & 1.7E-05 & 7000 & 18455.6094 & 2.1E-04 \\ \hline
 DGMRES & 324 & 166.7188 & 5.9E-08 & 608 & 1376.5000 & 6.7E-08 & 1202 & 10952.3438 & 6.6E-08 \\ \hline
\end{tabular}}
\end{table}


\begin{table}
\footnotesize \setlength{\abovecaptionskip}{0pt}
\setlength{\belowcaptionskip}{10pt} \centering{
\caption{\label{tab:all_wr_methods_table_dtMinus6_q2000_Lt25} The
number of iterations, computation time and the relative errors of
the approximate solutions for different methods on time interval
$[0,\ell_{t,J}\times J\times\Delta t]=[0,5\times 5\times\Delta t]$:
$\Delta t=10^{-6}$ and $q=2000$.}
\begin{tabular}{|l|ccc|ccc|ccc|}\hline
  & \multicolumn{3}{|c|}{$r=n^2,n=255$} & \multicolumn{3}{|c|}{$r=n^2,n=511$} & \multicolumn{3}{|c|}{$r=n^2,n=1023$} \\ \hline
 & IT & CPU & ERR & IT & CPU & ERR & IT & CPU & ERR \\ \hline
 WR-HSS & 35 & 63.7500 & 5.8E-11 & 17 & 150.1563 & 6.0E-11 & 14 & 487.7031 & 5.9E-10 \\ \hline
 WR-SOR & 98 & 18.8750 & 3.7E-10 & 1723 & 1202.9375 & 2.7E-10 & 1561 & 4166.5469 & 1.6E-10 \\ \hline
 DGMRES & 28 & 12.4531 & 7.8E-09 & 49 & 103.4063 & 6.0E-09 & 93 & 824.5313 & 3.7E-09 \\ \hline
\end{tabular}}
\end{table}


\begin{table}
\footnotesize \setlength{\abovecaptionskip}{0pt}
\setlength{\belowcaptionskip}{10pt} \centering{
\caption{\label{tab:all_wr_methods_table_dtMinus6_q3000_Lt25} The
number of iterations, computation time and the relative errors of
the approximate solutions for different methods on time interval
$[0,\ell_{t,J}\times J\times\Delta t]=[0,5\times 5\times\Delta t]$:
$\Delta t=10^{-6}$ and $q=3000$.}
\begin{tabular}{|l|ccc|ccc|ccc|}\hline
  & \multicolumn{3}{|c|}{$r=n^2,n=255$} & \multicolumn{3}{|c|}{$r=n^2,n=511$} & \multicolumn{3}{|c|}{$r=n^2,n=1023$} \\ \hline
 & IT & CPU & ERR & IT & CPU & ERR & IT & CPU & ERR \\ \hline
 WR-HSS & 23 & 41.0000 & 6.2E-11 & 13 & 114.1406 & 2.8E-10 & 15 & 547.4531 & 2.0E-10 \\ \hline
 WR-SOR & 427 & 81.5469 & 4.5E-10 & 1710 & 1199.2969 & 3.5E-10 & 6853 & 18363.9219 & 2.4E-10 \\ \hline
 DGMRES & 38 & 18.2500 & 2.6E-09 & 69 & 150.5313 & 3.1E-09 & 133 & 1204.5625 & 2.1E-09 \\ \hline
\end{tabular}}
\end{table}


\begin{table}
\footnotesize \setlength{\abovecaptionskip}{0pt}
\setlength{\belowcaptionskip}{10pt} \centering{
\caption{\label{tab:long_time_interval_dtMinus3} The number of
iterations and the relative errors of the approximate solutions for
WR-HSS method on time interval $[0,1]$: $\Delta t=10^{-3}$,
$\ell_t=\ell_{t,J} \times J=5\times 200$.}
\begin{tabular}{|l|cc|cc|}\hline
 & \multicolumn{2}{|c|}{$q=2000$} & \multicolumn{2}{|c|}{$q=3000$} \\ \hline
 & $r=n^2,n=127$ & $r=n^2,n=255$ & $r=n^2,n=127$ & $r=n^2,n=255$ \\ \hline
 IT &     97 &    116 &   116 &   133 \\ \hline
 ERR &    1.0E-06 &   7.8E-07 &   7.0E-07 &   7.8E-07  \\ \hline
\end{tabular}}
\end{table}


\begin{table}
\footnotesize \setlength{\abovecaptionskip}{0pt}
\setlength{\belowcaptionskip}{10pt} \centering{
\caption{\label{tab:long_time_interval_dtMinus4} The number of
iterations and the relative errors of the approximate solutions for
WR-HSS method on time interval $[0,1]$: $\Delta t=10^{-4}$,
$\ell_t=\ell_{t,J} \times J=5\times 2000$.}
\begin{tabular}{|l|cc|cc|}\hline
 & \multicolumn{2}{|c|}{$q=2000$} & \multicolumn{2}{|c|}{$q=3000$} \\ \hline
 & $r=n^2,n=127$ & $r=n^2,n=255$ & $r=n^2,n=127$ & $r=n^2,n=255$ \\ \hline
 IT &     86 &    106 &   108 &   126 \\ \hline
 ERR &    9.9E-08 &   6.3E-08 &   5.4E-08 &   5.8E-08 \\ \hline
\end{tabular}}
\end{table}


\begin{table}
\footnotesize \setlength{\abovecaptionskip}{0pt}
\setlength{\belowcaptionskip}{10pt} \centering{
\caption{\label{tab:long_time_interval_dtMinus5} The number of
iterations and the relative errors of the approximate solutions for
WR-HSS method on time interval $[0,1]$: $\Delta t=10^{-5}$,
$\ell_t=\ell_{t,J} \times J=5\times 20000$.}
\begin{tabular}{|l|cc|cc|}\hline
 & \multicolumn{2}{|c|}{$q=2000$} & \multicolumn{2}{|c|}{$q=3000$} \\ \hline
 & $r=n^2,n=127$ & $r=n^2,n=255$ & $r=n^2,n=127$ & $r=n^2,n=255$ \\ \hline
 IT &     25 &    30 &    38 &    45 \\ \hline
 ERR &    7.7E-09 &   6.3E-09 &   5.5E-09 &   4.8E-09 \\ \hline
\end{tabular}}
\end{table}



\begin{thebibliography}{99}

\bibitem{Bai09NLAA}
Z.Z. Bai,
\newblock
Optimal parameters in the HSS-like methods for
saddle-point problems,
\newblock
Numer. Linear Algebra Appl. 16(2009) 447-479.

\bibitem{BaiBenziChen10C}
Z.Z. Bai, M. Benzi, F. Chen,
\newblock
Modified HSS iteration methods for a class of complex symmetric
linear systems,
\newblock
Comput. 87(2010) 93-111.

\bibitem{BaiBenziChen11NA}
Z.Z. Bai, M. Benzi, F. Chen,
\newblock
On preconditioned MHSS iteration methods for complex symmetric
linear systems,
\newblock
Numer. Algorithms 56(2011) 297-317.

\bibitem{BaiBenziChenWang13IMA}
Z.Z. Bai, M. Benzi, F. Chen, Z.Q. Wang,
\newblock
Preconditioned MHSS iteration methods for a class of block
two-by-two linear systems with applications to distributed control
problems,
\newblock
IMA J. Numer. Anal. 33(2013) 343-369.

\bibitem{BaiGolub07IMA}
Z.Z. Bai, G.H. Golub,
\newblock
Accelerated Hermitian and skew-Hermitian splitting iteration methods
for saddle-point problems,
\newblock
IMA J. Numer. Anal. 27(2007) 1-23.

\bibitem{BaiLi06}
Z.Z. Bai, G.H. Golub, C.K. Li,
\newblock
Optimal parameter in Hermitian and skew-Hermitian splitting method
for certain two-by-two block matrices,
\newblock
SIAM J. Sci. Comput. 28(2006) 583-603.

\bibitem{BaiGolubLi07MC}
Z.Z. Bai, G.H. Golub, C.K. Li,
\newblock
Convergence properties of preconditioned Hermitian and
skew-Hermitian splitting methods for non-Hermitian positive
semidefinite matrices,
\newblock
Math. Comput. 76(2007) 287-298.

\bibitem{BaiYin05}
Z.Z. Bai, G.H. Golub, L.Z. Lu, J.F. Yin,
\newblock
Block triangular and skew-Hermitian splitting methods for
positive-definite linear systems,
\newblock
SIAM J. Sci. Comput. 26(2005) 844-863.

\bibitem{Bai-03}
Z.Z. Bai, G.H. Golub, M.K. Ng,
\newblock
Hermitian and skew-Hermitian splitting methods for non-Hermitian
positive definite linear systems,
\newblock
SIAM J. Matrix Anal. Appl. 24(2003) 603-626.

\bibitem{BaiGolubNg07NLAA}
Z.Z. Bai, G.H. Golub, M.K. Ng,
\newblock
On successive overrelaxation acceleration of the Hermitian and
skew-Hermitian splitting iterations,
\newblock
Numer. Linear Algebra Appl. 14(2007) 319-335.

\bibitem{BaiNg08}
Z.Z. Bai, G.H. Golub, M.K. Ng,
\newblock
On inexact Hermitian and skew-Hermitian splitting methods for
non-Hermitian positive definite linear systems,
\newblock
Linear Algebra Appl. 428(2008) 413-440.

\bibitem{BaiPan04}
Z.Z. Bai, G.H. Golub, J.Y. Pan,
\newblock
Preconditioned Hermitian and skew-Hermitian splitting methods for
non-Hermitian positive semidefinite linear systems,
\newblock
Numer. Math. 98(2004) 1-32.

\bibitem{wr-lcc-4}
Z.Z. Bai, M.K. Ng, J.Y. Pan,
\newblock
Alternating splitting waveform relaxation method
and its successive overrelaxation acceleration,
\newblock
Comput. Math. Appl. 49(2005) 157-170.

\bibitem{BertGolub05}
D. Bertaccini, G.H. Golub, S.S. Capizzano, C.T. Possio,
\newblock
Preconditioned HSS methods for the solution of non-Hermitian
positive definite linear systems and applications to the
discrete convection-diffusion equation,
\newblock
Numer. Math. 99(2005) 441-484.

\bibitem{SL-Camp-1}
S.L. Campbell,
\newblock
Singular Systems of Differential Equations,
\newblock
Pitman Advanced Publishing Program,
London and Melbourne, 1980.

\bibitem{Douglas-56}
J. Douglas Jr., H.H. Rachford Jr.,
\newblock
Alternating direction methods for three space variables,
\newblock
Numer. Math. 4(1956) 41-63.

\bibitem{Elman90}
H.C. Elman, G.H. Golub,
\newblock
Iterative methods for cyclically reduced non-self-adjoint linear systems,
\newblock
Math. Comput. 54(1990) 671-700.

\bibitem{Elman91}
H.C. Elman, G.H. Golub,
\newblock
Iterative methods for cyclically reduced non-self-adjoint linear systems II,
\newblock
Math. Comput. 56(1990) 215-242.

\bibitem{Elman-05}
H.C. Elman, D.J. Silvester, A.J. Wathen,
\newblock
Finite Elements and Fast Iterative Solvers: with Applications in
Incompressible Fluid Dynamics,
\newblock
Oxford University Press, New York, 2005.

\bibitem{Hackbusch85}
W. Hackbusch,
\newblock
Multi-Grid Mehtods and Applications,
\newblock
Springer-Verlag, Berlin Heidelberg, 1985.

\bibitem{Hochstadt-89}
H. Hochstadt,
\newblock
Integral Equations,
John Wiley \& Sons,
New York, 1989.

\bibitem{Jan-97}
J. Janssen, S. Vandewalle,
\newblock
On SOR waveform relaxation methods,
\newblock
SIAM J. Numer. Anal. 34(1997) 2456-2481.

\bibitem{Jiang-09}
Y.L. Jiang,
\newblock
Waveform Relaxation Methods,
\newblock
Science Press, Beijing, 2009.

\bibitem{LelRS-82}
E. Lelarasmee, A. Ruheli, A.L. Sangiovanni-Vincentelli,
\newblock
The waveform relaxation method for time-domain analysis of
large scale integrated circuits,
\newblock
IEEE Trans. Comput.-Aided Des. Integr. Circuits Syst.
CAD-1(1982) 131-145.

\bibitem{FL-Lewis}
F.L. Lewis,
\newblock
A survey of linear singular systems,
\newblock
Circuits Syst. Signal Process. 5(1986) 3-36.

\bibitem{wr-lcc-1}
U. Miekkala, O. Nevanlinna,
\newblock
Convergence of dynamic iteration methods for
initial value problems,
\newblock
SIAM J. Sci. Stat. Comput. 8(1987) 459-482.

\bibitem{wr-lcc-3}
J.Y. Pan, Z.Z. Bai,
\newblock
On the convergence of waveform relaxation
methods for linear initial value problems,
\newblock
J. Comput. Math. 22(2004) 681-698.

\bibitem{wr-lcc-5}
J.Y. Pan, Z.Z. Bai, M.K. Ng,
\newblock
Two-step waveform relaxation methods for
implicit linear initial value problems,
\newblock
Numer. Linear Algebra Appl. 12(2005) 293-304.

\bibitem{WangBai05}
J. Wang, Z.Z. Bai,
\newblock
Convergence analysis of two-stage waveform
relaxation method for the initial value problems,
\newblock
Appl. Math. Comput. 172(2006) 797-808.


\end{thebibliography}
\end{document}